\documentclass[12pt]{article}

\usepackage{geometry}
\usepackage{url}
\usepackage{comment}
\usepackage{indentfirst}  % 新增：让每个段落（包括首段）都缩进

\usepackage{amssymb,amsthm,amsmath,bm}
\usepackage{stmaryrd}

\usepackage{overpic}

\theoremstyle{remark}

\usepackage{graphicx}
\usepackage{epstopdf}
\usepackage{subfigure}
\usepackage{tabularx}
\usepackage{booktabs}
\usepackage{ragged2e}
\usepackage{siunitx}
\usepackage{caption}

\usepackage{listings}
\usepackage{xcolor}
\usepackage{algorithm}
\usepackage{algpseudocode}

\usepackage{tcolorbox}
\tcbuselibrary{breakable}

\newtcbox{\graycode}{
  on line,
  boxsep=0pt,
  left=1pt, right=1pt, top=1pt, bottom=1pt,
  colback=gray!10,
  colframe=gray!30,
  boxrule=0.5pt,
  arc=2pt,
  fontupper=\ttfamily
}

\usepackage[numbers,sort&compress]{natbib}
\usepackage{hyperref}

\providecommand{\email}[1]{\href{mailto:#1}{\texttt{#1}}}

\title{SaddleScape V1.0: A Python Package for Constructing Solution Landscapes via High-index Saddle Dynamics}
\author{
Yuyang Liu\thanks{School of Mathematical Sciences, Peking University, Beijing, 100871, China. ({liuyuyang@stu.pku.edu.cn}).}\and
Hua Su\thanks{Beijing International Center for Mathematical Research, Peking University, Beijing, 100871, China. ({suhua@pku.edu.cn}).}\and
Zixiang Xiao\thanks{School of Mathematical Sciences, Peking University, Beijing, 100871, China. ({xiaozixiang@stu.pku.edu.cn}).}\and
Lei Zhang\thanks{Beijing International Center for Mathematical Research, Center for Quantitative Biology, Center for Machine Learning Research, Peking University, Beijing, 100871, China. (\email{zhangl@math.pku.edu.cn}).}
\and Jin Zhao\thanks{Academy for Multidisciplinary Studies, Capital Normal University, and Beijing National Center for Applied Mathematics, Beijing, 100048, China.  (\email{zjin@cnu.edu.cn}).}
}

\date{\today}

\begin{document}

\maketitle

\begin{abstract}

We present \texttt{SaddleScape V1.0}, a Python software package designed for the exploration and construction of solution landscapes in complex systems. The package implements the High-index Saddle Dynamics (HiSD) framework and its variants, including the Generalized HiSD for non-gradient systems and the Accelerated HiSD. \texttt{SaddleScape V1.0} enables the systematic identification of critical points, including both local minima and high-index saddle points, by dynamically updating both the state estimate and an associated subspace characterizing the saddle's local manifold. It supports both gradient systems, defined by energy functions/functionals, and general non-gradient autonomous dynamical systems. Key features include automatic differentiation for symbolic inputs, numerical approximation techniques for Hessian-vector products, diverse eigenvalue solvers, and algorithms for constructing solution landscapes. The software offers a user-friendly interface with flexible parameter configuration, tools for trajectory and landscape visualization, and data export capabilities. By providing an efficient and accessible implementation of advanced saddle dynamics, \texttt{SaddleScape V1.0} facilitates the construction of solution landscapes, empowering researchers in various scientific disciplines to gain deeper insights into the hierarchical structure of complex systems. The source code is available at the repository \href{https://github.com/HiSDpackage/saddlescape}{https://github.com/HiSDpackage/saddlescape}. The package's introductory website is available at \href{https://hisdpackage.github.io/saddlescape}{https://hisdpackage.github.io/saddlescape}.

\end{abstract}

\hspace{-0.5cm}\textbf{Keywords:}
\small{SaddleScape, Solution Landscape, High-index Saddle Dynamics, Saddle Points}\\

\section{Introduction}

Many practical problems in physics, chemistry, computer science, and related fields can be formulated as minimization problems of multivariable nonlinear functions/functionals, representing amino acid positions in protein folding \cite{protein-folding1,protein-folding2,protein-folding3}, atomic coordinates in Lennard-Jones clusters \cite{atomic}, or parameters of artificial neural networks \cite{neural-networks1,neural-networks2,neural-networks3}. Each local minimum of the energy function/functional corresponds to a metastable state, and different minima are separated by energy barriers. The state located at the point of the minimal energy barrier between two minima is called the transition state. Transition states play crucial roles in many scientific fields, such as finding critical nuclei and transition pathways in phase transformations \cite{PhysRevLett.104.148301,PhysRevLett.98.265703}, computing transition rates in chemical reactions \cite{transition-rates1,transition-rates2}, and analyzing biological dynamics \cite{biological-applications1,biological-applications2,biological-applications3}.

Besides local minima and transition states, the stationary points of the energy function/functional
also include high-index saddle points. Mathematically, these stationary points are solutions of $\nabla E(\boldsymbol{x})=0$. For a nondegenerate saddle point, Morse Theory \cite{Milnor+1963} characterizes its nature by the Morse index, defined as the maximal dimension of a subspace on which the Hessian is negative definite. For example, a metastable state has Morse index $0$, and a transition state has Morse index $1$. In this paper, non-degenerate critical points with Morse index $k>0$ are referred to as index-$k$ saddle points (for convenience, minima are sometimes also called index-$0$ saddle points).

Thus, to gain a deeper understanding of the energy landscape, it is crucial not only to locate minima but also to identify the index-1 saddle points connecting them. Furthermore, lower-index and higher-index saddle points can be connected via dynamical search pathways produced by appropriate algorithms. Revealing these connections can provide profound insights into the structure of the energy landscape, thereby explaining various physical and chemical phenomena. However, most traditional numerical methods for searching saddle points, such as homotopy methods \cite{homotopy-methods1,homotopy-methods2} and deflation techniques \cite{deflation-techniques}, do not explicitly address the Morse index or the connectivity between saddle points. To address this, Zhang et al. proposed and analyzed the High-index Saddle Dynamics (HiSD) algorithm \cite{yin2019high}, suitable for searching high-index saddle points of a specified index. This algorithm enables systematic searches for saddle points of specified indices and reveals their interconnections.

Building on this foundation, Zhang et al. introduced the concept of solution landscape\cite{yin2020construction,yin2021searching}: a pathway map consisting of all the stationary points and their interconnections. This map describes a hierarchical structure in which saddle points are organized by index, showing how higher-index saddles relate to lower-index saddles. This hierarchical structure enhances our understanding of the relationships between saddle points of different indices and enables more systematic and comprehensive saddle point searches, thereby providing deeper insights into the energy landscape.

\begin{figure}[htbp]
    \centering
    \includegraphics[width=0.5\linewidth]{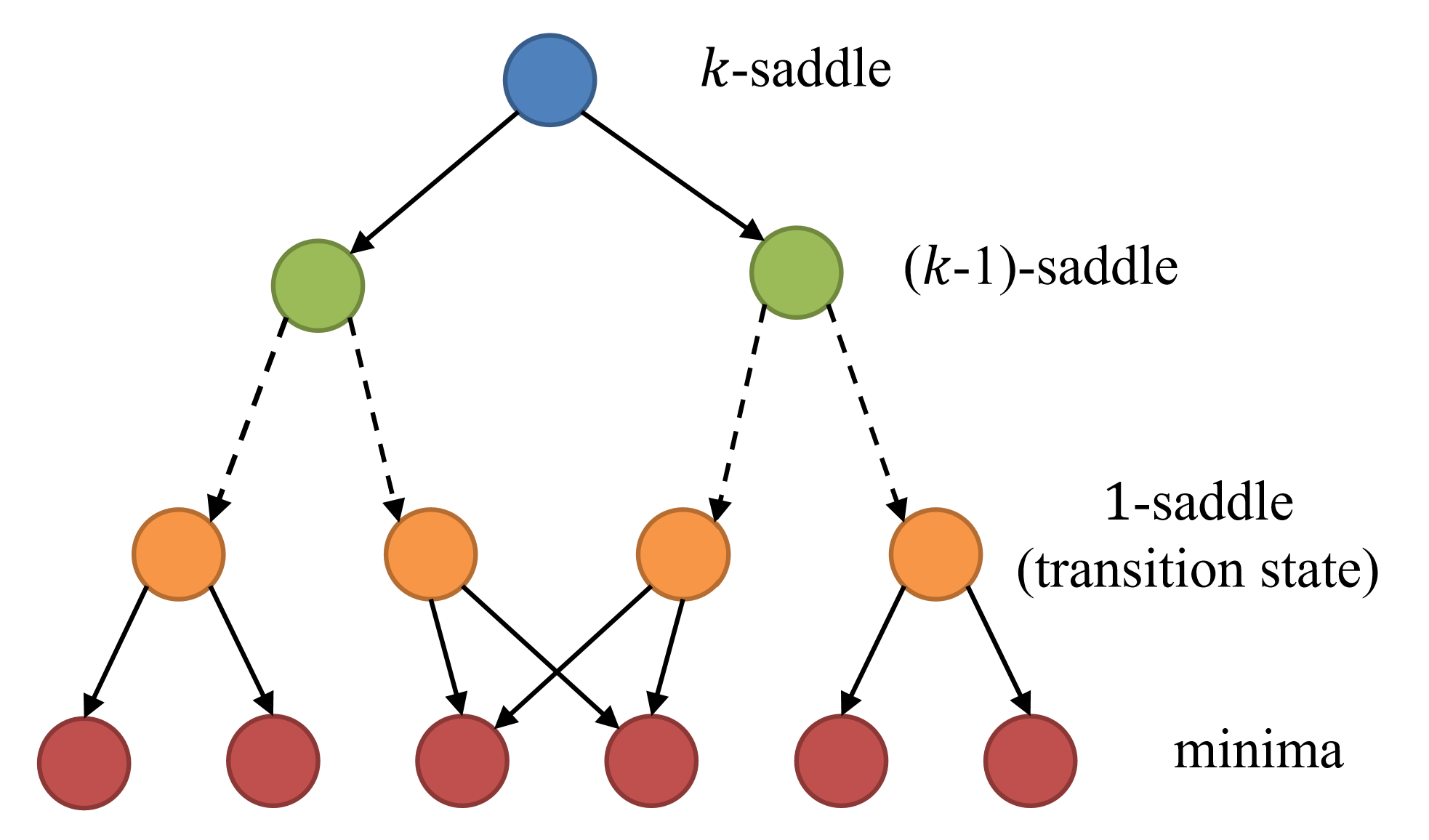}
    \caption{\textbf{A diagram of the solution landscape.} Nodes represent saddle points of different Morse indices, and directed edges indicate interconnections between saddles.}
    \label{fig:placeholder}
\end{figure}

Subsequent research extended the solution landscape concept and the HiSD method from energy landscapes (gradient systems) to non-gradient autonomous systems $\dot{\boldsymbol{x}}=\boldsymbol{F}(\boldsymbol{x})$, which are widely encountered in models describing physical, chemical, and biological dynamics, e.g., kinetic equations \cite{KineticEq}, reduced models related to the Navier–Stokes equations \cite{NavierStokes}, and biochemical reactions \cite{BioChem1,BioChem2,BioChem3}. This led to the proposal of Generalized High-index Saddle Dynamics (GHiSD) \cite{yin2021searching}. This framework generalizes the case of energy landscapes (gradient systems), where $\boldsymbol{F}(\boldsymbol{x})=-\nabla E(\boldsymbol{x})$. Furthermore, Constrained High-index Saddle Dynamics (CHiSD) \cite{yin2022constrained} was developed for systems with nonlinear equality constraints. These saddle dynamics algorithms and the construction methods of solution landscapes are systematically elaborated in the review article \cite{jcm2023}.

In recent years, the combination of HiSD with other numerical methods has provided more efficient implementations and solid theoretical support for constructing solution landscapes. In particular, representative variants include Accelerated High-index Saddle Dynamics (AHiSD) \cite{luo2025}, which introduces momentum to accelerate convergence; Improved High-index Saddle Dynamics (iHiSD) \cite{su2025sinum}, which combines saddle dynamics and gradient flow to provide theoretical justification for the global convergence of the algorithm; and data-driven HiSD variants that employ surrogate models constructed from sampled data, such as Gaussian-process-based models \cite{JJIAM2023} and neural-network-based models \cite{Liu2024Neural}, among others.  Moreover, researchers have combined different numerical schemes with HiSD to obtain various effective algorithms and conducted theoretical convergence and error analyses \cite{zhang2022sinum,luo2022sinum,csiam2023,miao2025,scm2023,luo2024NMPDE,cam2023,miao2025csiam,zhang2024probabilistic}.

Beyond the developments in theoretical and numerical methods, these algorithms have been successfully applied to study the structure of various systems. In condensed matter physics, the solution landscape concept and saddle dynamics has been used to describe crystal potential energy surfaces \cite{li2024acta}, to compute phase transition and nucleation of quasicrystals \cite{Yin2021transition,zhou2024nucleation}, and to construct and analyze defect landscapes in liquid crystal systems \cite{yin2020construction,han2021a,han2021solution,han2022elastic,yin2022solution,shi2022siap,shi2023nonlinearity,shi2024,shi2025siap,wang2021modeling}. In these complex systems, HiSD algorithms provide a systematic way to locate families of saddles and minima and to organize them into solution landscapes that clarify the interconnections of stationary states.

In theoretical physics, researchers have applied these algorithms to explore excited states of rotating Bose–Einstein condensates by mapping out the associated solution landscapes at different rotational frequencies \cite{yin2024revealing}. In biology, the solution landscape framework has been used to describe cell fate differentiation \cite{xue2024elife,Zhang2024.11.28.625944}, and to analyze cancer progression dynamics in tumor ecological microenvironments \cite{wang2025siap}. Related techniques have also been employed for other complex systems with multiple equilibria, including diblock copolymer-homopolymer models \cite{xu2021solution},  Cylinder-Gyroid interface systems described by Landau–Brazovskii models \cite{yao2022cicp}, space-fractional phase field models \cite{yu2022jcp}, multi-phase fluid flow systems \cite{jcp2023}, and droplet wetting and directional transport on rough surfaces \cite{zhang2025softmatter}. By constructing solution landscapes that reveal the connectivity between saddle points and minima in these settings, saddle dynamics methods enable systematic studies of saddle-point problems and provide detailed insight into the underlying physical mechanisms. These examples also highlight the need for practical computational frameworks that can construct such solution landscapes in a systematic and automated way.

Despite these methodological and application advances, the practical construction of solution landscapes in concrete problems often remains at the level of manually assembling individual connections between saddle points. In a typical workflow, one has to design and run separate HiSD-based searches for each edge in the landscape (i.e., each direct connection between two critical points), specify \emph{ad hoc} stopping criteria for different runs, and rely heavily on personal experience with HiSD and its variants in order to explore and refine the landscape. Moreover, different variants of HiSD (such as GHiSD, CHiSD, AHiSD, and data-driven extensions) typically adopt slightly different notations, conventions, and implementation details, which can easily lead to confusion and inconsistency in practice. At present, there is still a lack of a holistic, automated algorithmic framework for constructing solution landscapes, accompanied by a systematic, integrated, and user-friendly software implementation that brings these methods together in a single tool.

To facilitate research and application of solution landscapes, we have developed the \texttt{SaddleScape V1.0} software package, which implements a series of HiSD algorithms for constructing solution landscapes. The key contribution of \texttt{SaddleScape V1.0} is the implementation of a systematic, hierarchical, and automated procedure to construct solution landscapes. Our package also provides efficient and user-friendly implementation of HiSD algorithms, enabling researchers to easily explore and apply HiSD methods. Specifically, our software package offers the following key advantages:

\begin{itemize}
    \item \textbf{Automated Solution Landscape Construction:} We implement a systematic downward search procedure based on Breadth-First Search (BFS) with First-In First-Out (FIFO) queue scheduling, combined with randomized perturbations. Upward searches from minima or lower-index saddles are also supported through restart interfaces, which allow users to refine or extend an existing, partially constructed solution landscape.

    \item \textbf{Unified Framework for HiSD Variants and System Types:} The package provides a unified implementation of HiSD and its variants (such as GHiSD for non-gradient systems and AHiSD with momentum acceleration) within a single interface. It supports both gradient systems and general non-gradient autonomous systems, including functional systems after spatial discretization.

    \item \textbf{Flexible and Efficient Numerical Implementations:} Most steps of the solving process can be configured by the user. Our system accepts both numerical functions and symbolic expressions, supports custom gradients as well as automatic differentiation, and incorporates matrix-free Hessian/Jacobian-vector products, multiple eigenvalue solvers, and several step-size and acceleration strategies (including Barzilai–Borwein, heavy-ball, and Nesterov methods) to enhance robustness and efficiency.

    \item \textbf{Integrated Postprocessing and Visualization:}
    We provide wrapper methods for search trajectory visualization, solution landscape construction, and data saving in various formats. These postprocessing modules offer a quick and intuitive understanding of both the search process and the resulting solution landscape, and facilitate downstream analysis.

    \item \textbf{Modular Architecture and User-friendly Interface:}
    The package is organized into modular components for parameter checking, symbolic and numerical derivative handling, Hessian/Jacobian operators, eigenvalue solvers, iterative saddle dynamics solvers, and solution landscape construction. Users can configure system definitions, solver parameters, eigenvalue methods, and acceleration schemes via a clear set of options. Default values are provided for many parameters, lowering the barrier for new users while still allowing advanced users to fine-tune numerical details and experiment with extended algorithms by replacing the corresponding modules.
\end{itemize}

In summary, the package simplifies the saddle point search process, offers flexible parameter settings, and provides a range of visualization and data-handling tools that lower the barrier to constructing and analyzing solution landscapes in practice.

For ease of reading, definitions of commonly used symbols are listed in Table \ref{tab:notation}.
\begin{table}[htbp]
\centering
\caption{\textbf{Summary of mathematical symbols and key definitions used in the HiSD framework.}}
\label{tab:notation}
\begin{tabular}{|>{\centering\arraybackslash}m{2cm}|>{\raggedright\arraybackslash}m{\dimexpr\textwidth-3.5cm\relax}|}
\hline
\textbf{Symbol} & \multicolumn{1}{c|}{\textbf{Description}} \\
\hline
$\boldsymbol{a}$ & Bold font denotes vectors. \\
\hline
$\mathbb{A}$ & Blackboard bold denotes standard algebraic objects (e.g. matrices or number fields). In particular, $\mathbb{R}$ and $\mathbb{C}$ denote the real and complex number fields. \\
\hline
$\mathcal{A}$ & Calligraphic font denotes spaces. \\
\hline
$d$ & Dimension of the system. \\ 
\hline
$\langle \cdot,\cdot \rangle$ & Inner product in Hilbert space $\mathcal{H}$; written as $\boldsymbol{x}^{\top} \boldsymbol{y}$ in finite dimensions. \\ 
\hline
$E(\boldsymbol{x})$ & Twice Fréchet differentiable energy functional defined on a real Hilbert space $\mathcal{H}$. \\
\hline
$\boldsymbol{F}(\boldsymbol{x})$ & The vector field $\boldsymbol{F}(\boldsymbol{x})$; in gradient systems $\boldsymbol{F}(\boldsymbol{x})=-\nabla E(\boldsymbol{x})$. \\
\hline
$\mathbb{H}(\boldsymbol{x})$ & Hessian matrix: $\mathbb{H}(\boldsymbol{x}) = \nabla^2E(\boldsymbol{x})$. \\
\hline
$\mathbb{J}(\boldsymbol{x})$ & Jacobian matrix: $\mathbb{J}(\boldsymbol{x}) = \nabla \boldsymbol{F}(\boldsymbol{x})$. In gradient systems where $\boldsymbol{F}(\boldsymbol{x}) = -\nabla E(\boldsymbol{x})$, we have $\mathbb{J}(\boldsymbol{x}) = -\mathbb{H}(\boldsymbol{x})$. \\
\hline
$\boldsymbol{\hat{x}}$ & Critical point (gradient systems) or equilibrium point (non-gradient systems). \\
\hline
$\mathcal{W}^u(\hat{\boldsymbol{x}})$ & Unstable subspace spanned by eigenvectors with eigenvalues corresponding to positive real parts. \\
\hline
$\mathcal{W}^s(\hat{\boldsymbol{x}})$ & Stable subspace spanned by eigenvectors with eigenvalues corresponding to negative real parts. \\
\hline
$\mathcal{W}^c(\hat{\boldsymbol{x}})$ & Center subspace spanned by eigenvectors with eigenvalues corresponding to zero real parts. \\
\hline
$k_u$ & Dimension of unstable subspace (index of equilibrium in non-gradient systems). \\
\hline
$k_s$ & Dimension of stable subspace. \\
\hline
$k_c$ & Dimension of center subspace. \\
\hline
$k$-saddle & Index-$k$ saddle point (Morse index $k$ in gradient systems, unstable dimension $k_u = k$ in non-gradient systems). \\
\hline
\end{tabular}
\end{table}

\section{High-index Saddle Dynamics}\label{High-index Saddle Dynamics}
This section provides a brief introduction to a series of HiSD algorithms. For details, see our website \href{https://hisdpackage.github.io/saddlescape}{https://hisdpackage.github.io/saddlescape} or the corresponding references. These algorithms constitute the theoretical foundation of this software package. The content is divided into three parts:
\begin{itemize}
    \item The basic HiSD\cite{yin2019high} algorithm for gradient systems with twice Fréchet differentiable energy functions;
    \item The extension to non-gradient systems, resulting in the GHiSD\cite{yin2021searching} algorithm;
    \item The momentum-accelerated AHiSD\cite{luo2025} algorithm.
\end{itemize}

\subsection{HiSD for Gradient System}\label{HiSD: High-index Saddle Dynamics for Gradient System}  
\subsubsection{Description of the System}\label{Description of the System}
Given a twice Fr\'echet-differentiable energy functional $E(\boldsymbol{x})$ defined on a real Hilbert space $\mathcal{H}$ with inner product $\langle \cdot,\cdot \rangle$, the natural force is defined as $\boldsymbol{F}(\boldsymbol{x}) = -\nabla E(\boldsymbol{x})$ and the Hessian matrix as $\mathbb{H}(\boldsymbol{x}) = \nabla^2 E(\boldsymbol{x})$.

\begin{itemize}
  \item A point $\boldsymbol{\hat{x}} \in \mathcal{H}$ is a critical point of $E(\boldsymbol{x})$ if $\lVert \boldsymbol{F}(\boldsymbol{\hat{x}}) \rVert = 0$. A critical point is non-degenerate if its Hessian matrix $\mathbb{H}(\boldsymbol{\hat{x}})$ is invertible with a bounded inverse.
  \item A critical point that is not a local extremum is termed a saddle point. By convention, local minima are index-$0$ saddle points, and local maxima are index-$d$ saddle points, where $d$ denotes the system dimension.
  \item For a non-degenerate critical point $\boldsymbol{\hat{x}}$, its Morse index is the maximal dimension of a subspace $\mathcal{K} \subset \mathcal{H}$ such that $\mathbb{H}(\boldsymbol{\hat{x}})$ is negative definite on $\mathcal{K}$. The objective is to locate index-$k$ saddle points ($k$-saddles) on the potential energy surface (PES).
  \item For simplicity, we assume that the dimension of $\mathcal{H}$ is $d$ and write the inner product $\langle \boldsymbol{x}, \boldsymbol{y} \rangle$ as $\boldsymbol{x}^{\top} \boldsymbol{y}$.
\end{itemize}

Equivalently, the index of a saddle point corresponds to the number of negative eigenvalues of its Hessian matrix. In the context of congruence canonical forms for real symmetric matrices, this quantity is also known as the negative index of inertia (the number of negative eigenvalues).

\subsubsection{Derivation of Index-$k$ HiSD}

For a non-degenerate critical point $\boldsymbol{\hat{x}}$, let $\hat{\mathcal{V}}$ be the subspace spanned by the eigenvectors of $\mathbb{H}(\boldsymbol{\hat{x}})$ corresponding to its negative eigenvalues. Since $\mathbb{H}(\boldsymbol{\hat{x}})$ is negative definite on $\hat{\mathcal{V}}$ and positive definite on $\hat{\mathcal{V}}^{\perp}$, $\boldsymbol{\hat{x}}$ is a local maximum in the affine space $\boldsymbol{\hat{x}} + \hat{\mathcal{V}}$ and a local minimum in $\boldsymbol{\hat{x}} + \hat{\mathcal{V}}^{\perp}$.

The projections of $\boldsymbol{\hat{x}}$ onto $\hat{\mathcal{V}}$ and $\hat{\mathcal{V}}^{\perp}$ are denoted by $\boldsymbol{\hat{x}}_{\hat{\mathcal{V}}}$ and $\boldsymbol{\hat{x}}_{\hat{\mathcal{V}}^{\perp}}$, respectively. The pair $(\boldsymbol{v}, \boldsymbol{w}) = (\boldsymbol{\hat{x}}_{\hat{\mathcal{V}}}, \boldsymbol{\hat{x}}_{\hat{\mathcal{V}}^{\perp}})$ satisfies the minimax problem
$$
\min_{\boldsymbol{w} \in \hat{\mathcal{V}}^{\perp}} \max_{\boldsymbol{v} \in \hat{\mathcal{V}}} E(\boldsymbol{v} + \boldsymbol{w}).
$$
As $\hat{\mathcal{V}}$ is unknown, one obtains an iterative scheme that updates both the current iterate $\boldsymbol{x}$ (an estimate of the critical point) and an approximation $\mathcal{V}$ to $\hat{\mathcal{V}}$.

\noindent\textbf{Dynamics of $\boldsymbol{x}$}

Intuitively, the update of $\boldsymbol{x}$ is designed such that the projection of $\boldsymbol{\dot{x}}$ onto the space $\mathcal{V}$, $\mathcal{P}_{\mathcal{V}}\boldsymbol{\dot{x}}$, is the ascent direction of the energy function $E(\boldsymbol{x})$, while the projection onto the complementary space $\mathcal{V}^{\perp}$, $\mathcal{P}_{\mathcal{V}^{\perp}}\boldsymbol{\dot{x}}$, is the descent direction. Here, $\mathcal{P}_{\mathcal{V}}$ and $\mathcal{P}_{\mathcal{V}^{\perp}}$ denote the orthogonal projection operators onto the subspace $\mathcal{V}$ and its orthogonal complement $\mathcal{V}^{\perp}$, respectively.

Noting that the natural force $\boldsymbol{F}(\boldsymbol{x})=-\nabla E(\boldsymbol{x})$ is the steepest descent direction, one may set
$$
\mathcal{P}_{\mathcal{V}}\boldsymbol{\dot{x}}=-\mathcal{P}_{\mathcal{V}}\boldsymbol{F}(\boldsymbol{x})
,\quad
\mathcal{P}_{\mathcal{V}^{\perp}}\boldsymbol{\dot{x}}=\mathcal{P}_{\mathcal{V}^{\perp}}\boldsymbol{F}(\boldsymbol{x})=\boldsymbol{F}(\boldsymbol{x})-\mathcal{P}_{\mathcal{V}}\boldsymbol{F}(\boldsymbol{x}).
$$
Then, with two positive relaxation parameters $\beta_{\mathcal{V}}$ and $\beta_{\mathcal{V}^{\perp}}$, the dynamics of $\boldsymbol{x}$ can be given as
\[
\boldsymbol{\dot{x}}=\beta_{\mathcal{V}}(-\mathcal{P}_{\mathcal{V}}\boldsymbol{F}(\boldsymbol{x}))+\beta_{\mathcal{V}^{\perp}}(\boldsymbol{F}(\boldsymbol{x})-\mathcal{P}_{\mathcal{V}}\boldsymbol{F}(\boldsymbol{x})).
\]
For simplicity, we set $\beta_{\mathcal{V}}=\beta_{\mathcal{V}^{\perp}}=\beta>0$ and take an orthonormal basis $\{\boldsymbol{v}_i\}_{i=1}^k$ for $\mathcal{V}$ with each $\boldsymbol{v}_i$ taken as a column vector, the above equation reduces to
\begin{equation}\label{2.1}
\beta^{-1}\boldsymbol{\dot{x}}=\boldsymbol{F}(\boldsymbol{x})-2\mathcal{P}_{\mathcal{V}}\boldsymbol{F}(\boldsymbol{x})=\left( \mathbb{I} - 2 \sum_{i=1}^{k} \boldsymbol{v}_i \boldsymbol{v}_i^\top \right) \boldsymbol{F}(\boldsymbol{x}).
\end{equation}

\noindent\textbf{Update of Subspace $\mathcal{V}$}

To approximate $\hat{\mathcal{V}}$, define $\mathcal{V}$ as the span of eigenvectors corresponding to the $k$ smallest eigenvalues of $\mathbb{H}(\boldsymbol{x})$. This subspace is computed via the simultaneous Rayleigh-quotient iterative minimization (SIRQIT) method. Consider solving the following $k$ constrained optimization problems simultaneously ($i=1,2,\ldots,k$):
\begin{equation}\label{Rayleigh Quotient Optimization}
\min_{\boldsymbol{v}_i} \langle \boldsymbol{v}_i,\mathbb{H}(\boldsymbol{x})\boldsymbol{v}_i \rangle, \hspace{4em} \text{s.t.} \hspace{1em} \langle \boldsymbol{v}_i,\boldsymbol{v}_j \rangle=\delta_{ij} \hspace{1em}
j=1,2,\ldots,i,
\end{equation}
where
\[
\delta_{ij}=
\begin{cases}
    1 & \text{if } i=j, \\
    0 & \text{if } i \neq j .
\end{cases}
\]
According to the method of Lagrange multipliers, taking the negative gradient of the Rayleigh quotient as the driving term and enforcing the orthonormality constraints via Lagrange multipliers yields the evolution equation of $\boldsymbol{\dot{v}}_i$ with Lagrange multipliers $\xi_j^{(i)}$ and relaxation parameter $\gamma>0$:
\begin{equation}\label{the dynamics of v_i with undetermined coefficient}
\gamma^{-1}\boldsymbol{\dot{v}}_i
 =-\mathbb{H}(\boldsymbol{x})\boldsymbol{v}_i+\displaystyle \sum_{j=1}^{i}\xi_j^{(i)}\boldsymbol{v}_j .
\end{equation}
Since the constraints $\left\langle \boldsymbol{v}_i, \boldsymbol{v}_j \right\rangle = \delta_{ij}$ hold at any time, we have
\begin{equation}\label{Constraint}
\left\langle \boldsymbol{\dot{v}}_i, \boldsymbol{v}_j \right\rangle + \left\langle \boldsymbol{v}_i, \boldsymbol{\dot{v}}_j \right\rangle = 0 ,
\quad i, j = 1, \ldots, k.
\end{equation}
Substituting Eq. \eqref{the dynamics of v_i with undetermined coefficient} into Eq. \eqref{Constraint} and solving for the Lagrange multipliers in $\boldsymbol{\dot{v}}_i$ sequentially for $i=1,2,\ldots,k$, we obtain
\[
\xi^{(i)}_i=\langle \boldsymbol{v}_i,\mathbb{H}(\boldsymbol{x})\boldsymbol{v}_i \rangle
,\quad
\xi^{(i)}_j=2\langle \boldsymbol{v}_j,\mathbb{H}(\boldsymbol{x})\boldsymbol{v}_i \rangle
\hspace{1em} j=1,2,\ldots,i-1.
\]  
Substituting back into Eq. \eqref{the dynamics of v_i with undetermined coefficient} yields the final dynamics for $\boldsymbol{v}_i$:
\begin{equation}\label{2.3}
  \gamma^{-1} \boldsymbol{\dot{v}}_i = -\left( \mathbb{I} - \boldsymbol{v}_i \boldsymbol{v}_i^\top - 2 \sum_{j=1}^{i-1} \boldsymbol{v}_j \boldsymbol{v}_j^\top \right) \mathbb{H}(\boldsymbol{x}) \boldsymbol{v}_i, \quad i = 1, \ldots, k.
\end{equation}

\noindent\textbf{Overall Dynamical System}

The complete HiSD dynamical system combines \eqref{2.1} and \eqref{2.3}:
\begin{equation}\label{2.4}
\begin{cases}
  \beta^{-1} \boldsymbol{\dot{x}} = \left( \mathbb{I} - 2 \sum_{i=1}^{k} \boldsymbol{v}_i \boldsymbol{v}_i^\top \right) \boldsymbol{F}(\boldsymbol{x}), \\
\gamma^{-1} \boldsymbol{\dot{v}}_i = -\left( \mathbb{I} - \boldsymbol{v}_i \boldsymbol{v}_i^\top - 2 \sum_{j=1}^{i-1} \boldsymbol{v}_j \boldsymbol{v}_j^\top \right) \mathbb{H}(\boldsymbol{x}) \boldsymbol{v}_i, \quad i = 1, \ldots, k.
\end{cases}
\end{equation}
This framework dynamically updates both the critical point estimate and the associated subspace to locate saddle points. The linear stability analysis of these dynamics and further details are given in \cite{yin2019high}.

\subsubsection{Discretization of HiSD Method}\label{Discretization of HiSD method}
\noindent\textbf{SIRQIT with Approximate Hessian}

When the exact Hessian is available, Eq.\eqref{2.4} can be directly discretized for numerical solution. For many problems, however, computing the Hessian matrix is infeasible or prohibitively expensive, necessitating numerical approximation methods. Crucially, Eq.\eqref{2.4} requires only Hessian-vector products $\mathbb{H}(\boldsymbol{x}) \boldsymbol{v}_i$, which can be efficiently computed via finite differences.

The Hessian-vector product $\mathbb{H}(\boldsymbol{x}) \boldsymbol{v}$ can be approximated using the Taylor expansions:
$$
-\boldsymbol{F}(\boldsymbol{x}\pm l\boldsymbol{v}) = \nabla E(\boldsymbol{x}\pm l\boldsymbol{v}) = \nabla E(\boldsymbol{x}) \pm \mathbb{H}(\boldsymbol{x}) l \boldsymbol{v} + \mathcal{O}(\lVert l \boldsymbol{v} \rVert^2) .
$$
Combining these two expansions yields the approximation (also known as the dimer method), we get
\begin{equation}\label{the approximation of Hessian}
\boldsymbol{H}(\boldsymbol{x}, \boldsymbol{v}, l) = -\frac{\boldsymbol{F}(\boldsymbol{x}+l\boldsymbol{v}) - \boldsymbol{F}(\boldsymbol{x}-l\boldsymbol{v})}{2l} \approx \mathbb{H}(\boldsymbol{x}) \boldsymbol{v}.
\end{equation}
When $l\to 0$, theoretically exact gradients are obtained, but in practical numerical experiments, excessively small $l$ values should be avoided to prevent numerical instability. Incorporating this, the full system with approximate Hessian is obtained:
\begin{equation}
\begin{cases}
  \beta^{-1} \dot{\boldsymbol{x}} = \left( \mathbb{I} - 2 \sum_{i=1}^{k} \boldsymbol{v}_i \boldsymbol{v}_i^\top \right) \boldsymbol{F}(\boldsymbol{x}), \\
  \gamma^{-1} \dot{\boldsymbol{v}}_i = - \left( \mathbb{I} - \boldsymbol{v}_i \boldsymbol{v}_i^\top - 2 \sum_{j=1}^{i-1} \boldsymbol{v}_j \boldsymbol{v}_j^\top \right) \boldsymbol{H}(\boldsymbol{x}, \boldsymbol{v}_i, l) \quad i=1,\ldots,k .
\end{cases}
\end{equation}
Here, $l$ is taken as a fixed constant, and its selection can be determined by considering both the target accuracy and numerical stability.

Direct discretization of the above dynamics gives the following algorithm:
\begin{algorithm}[H]
\caption{\textbf{HiSD with approximate Hessian for a $k$-saddle.}}
\label{alg:hisd}

\textbf{Input :} $k \in \mathbb{N},\ l > 0,\ \boldsymbol{x}^{(0)} \in \mathcal{H},\ \{\boldsymbol{v}_i^{(0)}\}_{i=1}^k \subset \mathcal{H}$ satisfying $\left\langle \boldsymbol{v}_j^{(0)}, \boldsymbol{v}_i^{(0)} \right\rangle = \delta_{ij}$.

\begin{algorithmic}[1]
\State Set $n=0$, compute $\boldsymbol{f}^{(0)} = \boldsymbol{F}(\boldsymbol{x}^{(0)})$;
\Repeat
  \State $\boldsymbol{x}^{(n+1)} = \boldsymbol{x}^{(n)} + \beta^{(n)} \boldsymbol{g}^{(n)}$;
  \For{$i = 1:k$} \do \\
  \State $\boldsymbol{v}_i^* = \boldsymbol{v}_i^{(n)} + \gamma_i^{(n)} \boldsymbol{d}_i^{(n)}$;
  \State $\boldsymbol{v}_i^* = \boldsymbol{v}_i^* - \sum_{j=1}^{i-1} \left\langle \boldsymbol{v}_j^{(n+1)}, \boldsymbol{v}_i^* \right\rangle \boldsymbol{v}_j^{(n+1)}$;
  \State $\boldsymbol{v}_i^{(n+1)} = \boldsymbol{v}_i^* / \|\boldsymbol{v}_i^*\|$;
  \EndFor
  \State $\boldsymbol{f}^{(n+1)} = \boldsymbol{F}(\boldsymbol{x}^{(n+1)})$;
  \State $n := n + 1$;
\Until{$\| \boldsymbol{f}^{(n)} \| < \varepsilon_{\boldsymbol{F}}$}
\State \textbf{Output :} $\boldsymbol{x}^{(n)}, \boldsymbol{v}_1^{(n)}, \ldots, \boldsymbol{v}_k^{(n)}$.
\end{algorithmic}
\end{algorithm}
\noindent where
$$
\boldsymbol{g}^{(n)} = \boldsymbol{f}^{(n)} - 2 \sum_{i=1}^k \left\langle \boldsymbol{v}_i^{(n)}, \boldsymbol{f}^{(n)} \right\rangle \boldsymbol{v}_i^{(n)},
$$

$$
\boldsymbol{u}_i^{(n)} = \boldsymbol{H}\left( \boldsymbol{x}^{(n+1)}, \boldsymbol{v}_i^{(n)}, l \right),
$$

$$
\boldsymbol{d}_i^{(n)} = -\boldsymbol{u}_i^{(n)} + \left\langle \boldsymbol{v}_i^{(n)}, \boldsymbol{u}_i^{(n)} \right\rangle \boldsymbol{v}_i^{(n)} + \sum_{j=1}^{i-1} 2 \left\langle \boldsymbol{v}_j^{(n)}, \boldsymbol{u}_i^{(n)} \right\rangle \boldsymbol{v}_j^{(n)}.
$$

The selection methods for the step sizes $\beta^{(n)}$ and $\gamma_i^{(n)}$ are discussed in Section \ref{Selection of Iterative Step Sizes}. Additionally, although the continuous dynamics preserve the orthonormality of $\boldsymbol{v}_i$, due to discretization with a finite time step, it is recommended to apply Gram–Schmidt orthonormalization in lines 6 and 7 of the above algorithm to enforce the orthonormal constraints. In practice, when updating $\boldsymbol{v}_i^{(n)}$ one can reduce the time step and increase the number of substeps; i.e., repeat lines 4–8 of the algorithm $N$ times with step size $\dfrac{\gamma_i^{(n)}}{N}$, which yields a finer update of the eigensubspace.

\noindent\textbf{Combining with LOBPCG Method}

Alternatively, $\mathcal{V}$ can be updated by using established eigensolvers such as the Locally Optimal Block Preconditioned Conjugate Gradient (LOBPCG) method, which efficiently computes the $k$ smallest eigenvectors with an attractive accuracy-efficiency trade-off.

The core idea of the algorithm is to approximate the Rayleigh quotient optimization \eqref{Rayleigh Quotient Optimization} on a subspace and gradually enrich this subspace to improve the approximation, ultimately converging to the true $k$ smallest eigenvectors. This subspace reduction transforms the original problem into a lower-dimensional, more tractable eigenvalue problem.

For details on the implementation, we refer the reader to \cite{yin2019high}. Here, we only note that in the numerical implementation of LOBPCG, the multiplication of the Hessian matrix by a matrix can be decomposed into multiple Hessian matrix-vector products. These products can then be approximated using Eq. \eqref{the approximation of Hessian} to enhance computational efficiency.

\subsubsection{Selection of Iterative Step Sizes}\label{Selection of Iterative Step Sizes}
Briefly, typical methods for selecting step sizes are as follows:

\noindent\textbf{Explicit Euler Method}

Set $\beta^{(n)} = \beta \Delta t,\gamma_i^{(n)} = \gamma \Delta t$ as fixed constants with $\beta,\gamma,\Delta t>0$. Here $\Delta t$ should not be too small to avoid slow convergence, nor too large to prevent algorithmic divergence.

\noindent\textbf{BB Gradient Method}

For the selection of $\beta^{(n)}$, we employ the following Barzilai–Borwein (BB) step sizes\cite{Barzilai1988}:

Let $\Delta \boldsymbol{x}^{(n)} = \boldsymbol{x}^{(n)} - \boldsymbol{x}^{(n-1)}$ and $\Delta \boldsymbol{g}^{(n)} = \boldsymbol{g}^{(n)} - \boldsymbol{g}^{(n-1)}$. Following the quasi-Newton framework, we minimize the following norms
$$
\min_{\beta^{(n)}} \|\Delta \boldsymbol{x}^{(n)} / \beta^{(n)} - \boldsymbol{g}^{(n)}\|,
\quad
\min_{\beta^{(n)}} \|\Delta \boldsymbol{x}^{(n)} - \beta^{(n)} \boldsymbol{g}^{(n)}\|.
$$
to obtain the BB1 and BB2 step sizes:
$$
\beta_{\text{BB1}}^{(n)} = \frac{\langle \Delta \boldsymbol{x}^{(n)}, \Delta \boldsymbol{x}^{(n)} \rangle}{\langle \Delta \boldsymbol{x}^{(n)}, \Delta \boldsymbol{g}^{(n)} \rangle},
\quad
\beta_{\text{BB2}}^{(n)} = \frac{\langle \Delta \boldsymbol{x}^{(n)}, \Delta \boldsymbol{g}^{(n)} \rangle}{\langle \Delta \boldsymbol{g}^{(n)}, \Delta \boldsymbol{g}^{(n)} \rangle}.
$$
In this context, the denominator of the BB1 step size may become very small or even zero, leading to numerical instability; hence, we adopt the BB2 step size. Additionally, we can set an upper bound $\tau$ on $\beta^{(n)} \|\boldsymbol{g}^{(n)}\|$ and take the absolute value of $\beta_{\text{BB2}}^{(n)}$ to avoid negative step sizes. Combining these constraints yields:
\begin{equation}\label{BB2step}
\beta^{(n)} = \min \left\{ \frac{\tau}{\|\boldsymbol{g}^{(n)}\|}, \left| \frac{\langle \Delta \boldsymbol{x}^{(n)}, \Delta \boldsymbol{g}^{(n)} \rangle}{\langle \Delta \boldsymbol{g}^{(n)}, \Delta \boldsymbol{g}^{(n)} \rangle} \right| \right\}.    
\end{equation}

\subsection{GHiSD for Non-Gradient System}\label{GHiSD: Generalized High-index Saddle Dynamics for Non-Gradient System}

\subsubsection{Description of $k$-saddle in Non-Gradient Autonomous Dynamical Systems}\label{Description of k-saddle in Non-Gradient Autonomous Dynamical Systems}
The HiSD algorithm operates on a twice Fr\'echet differentiable energy functional $E(\boldsymbol{x})$. In contrast, the GHiSD algorithm targets index-$k$ saddle points in general $d$-dimensional non-gradient autonomous dynamical systems:
\begin{equation}
  \boldsymbol{\dot{x}} = \boldsymbol{F}(\boldsymbol{x}),\quad \boldsymbol{x} \in \mathbb{R}^d,\quad \boldsymbol{F} \in \mathcal{C}^r(\mathbb{R}^d,\mathbb{R}^d),\quad r \geq 2.
\label{dynamical system}
\end{equation}

We first define the index of a saddle point in a non-gradient system.

In the context of dynamical systems, denote the Jacobian matrix of $\boldsymbol{F}(\boldsymbol{x})$ as $\mathbb{J}(\boldsymbol{x}) = \nabla \boldsymbol{F}(\boldsymbol{x})$, with the inner product $\langle \boldsymbol{x},\boldsymbol{y} \rangle = \boldsymbol{x}^{\top}\boldsymbol{y}$. A point $\boldsymbol{\hat{x}}$ satisfying $\boldsymbol{F}(\boldsymbol{\hat{x}}) = \boldsymbol{0}$ is termed an equilibrium (or stationary) point.

Linearizing the system near equilibrium by setting $\boldsymbol{x} = \boldsymbol{\hat{x}} + \boldsymbol{y}$ and substituting into \eqref{dynamical system} yields
$$
\boldsymbol{\dot{y}} = \mathbb{J}(\boldsymbol{\hat{x}})\boldsymbol{y} + \mathcal{O}(\|\boldsymbol{y}\|^2).
$$
Neglecting higher-order terms gives the linearized system:
$$
\boldsymbol{\dot{y}} = \mathbb{J}(\boldsymbol{\hat{x}})\boldsymbol{y}.
$$
We classify the right (generalized) eigenvectors of $\mathbb{J}(\hat{\boldsymbol{x}})$ by the sign of their eigenvalues' real parts, which defines three key subspaces:
\begin{itemize}
  \item \textbf{Unstable subspace $\mathcal{W}^u(\hat{\boldsymbol{x}})$:} Spanned by eigenvectors with eigenvalues having positive real parts, whose dimension is denoted by $k_u$. Perturbations along these directions grow exponentially, driving the system away from equilibrium.
  
  \item \textbf{Stable subspace $\mathcal{W}^s(\hat{\boldsymbol{x}})$:} Spanned by eigenvectors with eigenvalues of negative real part, whose dimension is denoted by $k_s$. Perturbations here decay exponentially toward equilibrium.
  
  \item \textbf{Center subspace $\mathcal{W}^c(\hat{\boldsymbol{x}})$:} Spanned by eigenvectors with eigenvalues of zero real part, whose dimension is denoted by $k_c$. Perturbations exhibit neutral dynamics, such as oscillations, without exponential growth or decay.
\end{itemize}

By the spectral decomposition of $\mathbb{J}(\hat{\boldsymbol{x}})$, $\mathbb{R}^d = \mathcal{W}^u \oplus \mathcal{W}^s \oplus \mathcal{W}^c$, with $k_u + k_s + k_c = d$. When all eigenvalues have non-zero real parts ($\mathbb{R}^d = \mathcal{W}^u \oplus \mathcal{W}^s$), the equilibrium is hyperbolic, which is our primary focus.

In gradient systems, where $\boldsymbol{F} = -\nabla E$, the relation $\mathbb{J}=-\mathbb{H}$ connects the Hessian to the Jacobian. The index $k_u$ of an equilibrium, defined as the dimension of its unstable subspace, generalizes the concept of the Morse index in gradient systems. Thus, we refer to index-$k_u$ equilibria as index-$k$ saddle points, with sinks ($k_u=0$) and sources ($k_u=d$) as special cases. Subsequent sections detail the search for these $k$-saddles in non-gradient systems.

\subsubsection{Derivation of Index-$k$ GHiSD}
\noindent\textbf{Dynamics of $\boldsymbol{x}$}

The construction principle for the $\boldsymbol{x}$-dynamics parallels that of HiSD. Near an equilibrium point, following the system force drives convergence toward the equilibrium along the stable subspace and divergence away from it along the unstable subspace. Therefore, we preserve the force direction in the stable subspace and reverse its direction in the unstable subspace. Since the exact unstable subspace $\mathcal{W}^{\mathrm{u}}(\boldsymbol{\hat{x}})$ at the $k$-saddle is unknown, analogous to HiSD, we define $\mathcal{V}$ to approximate $\mathcal{W}^{\mathrm{u}}(\boldsymbol{\hat{x}})$ as the real invariant subspace spanned by the $k$ (generalized) eigenvectors of $\mathbb
{J}(\boldsymbol{x})$ associated with the $k$ eigenvalues having the largest real parts. If these eigenvalues (and eigenvectors) are complex, $\mathcal{V}$ is taken as the real span obtained from the real and imaginary parts of the corresponding complex generalized eigenvectors. This yields
\begin{equation}\label{GHiSD the dynamics of x}
\beta^{-1}\boldsymbol{\dot{x}} = \left( \mathbb{I} - 2\mathcal{P}_{\mathcal{V}} \right) \boldsymbol{F}(\boldsymbol{x})=\left( \mathbb{I} - 2 \sum_{j=1}^{k} \boldsymbol{v}_j \boldsymbol{v}_j^\top \right) \boldsymbol{F}(\boldsymbol{x}).
\end{equation}
This has the same form as \eqref{2.1}, where $\{\boldsymbol{v}_i\}_{i=1}^k$ is an orthonormal basis for $\mathcal{V}$, and $\beta>0$ is a relaxation parameter.

\noindent\textbf{Update of Subspace $\mathcal{V}$}

While efficient Rayleigh quotient optimization is not applicable to non-symmetric matrices, note that the dynamics $\boldsymbol{\dot{v}} = \mathbb{J}\boldsymbol{v}$ naturally align $\boldsymbol{v}$ with the invariant subspace of $\mathbb{J}$ corresponding to eigenvalues with largest real parts. Therefore, one can similarly propose a dynamical equation with undetermined coefficients $\xi_j^{(i)}$ and relaxation constant $\gamma>0$:
$$
\gamma^{-1}\boldsymbol{\dot{v}}_i = \mathbb{J}(\boldsymbol{x}) \boldsymbol{v}_i + \sum_{j=1}^{i} \xi^{(i)}_j \boldsymbol{v}_j.
$$
Similar to HiSD, using the orthonormality constraints \eqref{Constraint}, we obtain
$$
\xi_{i}^{(i)} = -\langle \mathbb{J}(\boldsymbol{x}) \boldsymbol{v}_i, \boldsymbol{v}_i \rangle,
\quad\xi_{j}^{(i)} = -\langle \mathbb{J}(\boldsymbol{x}) \boldsymbol{v}_i, \boldsymbol{v}_j \rangle - \langle \boldsymbol{v}_i, \mathbb{J}(\boldsymbol{x}) \boldsymbol{v}_j \rangle, \quad j = 1, \cdots, i-1.
$$
Combining these results with Eq. \eqref{GHiSD the dynamics of x} yields the complete dynamical system:
\begin{equation}\label{Entire GHiSD}
\begin{cases} 
\beta^{-1}\dot{\boldsymbol{x}} = \left( \mathbb{I} - 2 \displaystyle \sum_{j=1}^{k} \boldsymbol{v}_j \boldsymbol{v}_j^\top \right) \boldsymbol{F}(\boldsymbol{x}),\\
\gamma^{-1}\dot{\boldsymbol{v}}_i = \left( \mathbb{I} - \boldsymbol{v}_i \boldsymbol{v}_i^\top \right) \mathbb{J}(\boldsymbol{x}) \boldsymbol{v}_i - \displaystyle \sum_{j=1}^{i-1} \boldsymbol{v}_j \boldsymbol{v}_j^\top \left( \mathbb{J}(\boldsymbol{x}) + \mathbb{J}^\top(\boldsymbol{x}) \right) \boldsymbol{v}_i, \quad i = 1, \cdots, k.
\end{cases}
\end{equation}
The linear stability of this system is discussed in \cite{yin2021searching}.

\subsubsection{A Simplified Discrete Form of GHiSD}\label{Dis_GHiSD}

For the algorithmic implementation of the GHiSD method, similar to the HiSD algorithm, explicit Euler discretization of \eqref{Entire GHiSD} requires additional enforced orthonormalization operations due to discretization errors caused by the time step. Subsequent research \cite{zhang2025cms} has shown that the explicit Euler discretization of \eqref{Entire GHiSD} is equivalent to the following more concise discretization method, which is favored in practical computations due to its higher efficiency:
\begin{equation}\label{discrete GHiSD of W^u}
\begin{cases}
\tilde{\boldsymbol{v}}_i^{(n+1)} = \boldsymbol{v}_i^{(n)} + \gamma \mathbb{J}(\boldsymbol{x}) \boldsymbol{v}_i^{(n)} \quad i = 1, \ldots, k, \\
\left[ \boldsymbol{v}_1^{(n+1)}, \ldots, \boldsymbol{v}_k^{(n+1)} \right] = \text{orth} \left( \tilde{\boldsymbol{v}}_1^{(n+1)}, \ldots, \tilde{\boldsymbol{v}}_k^{(n+1)} \right) ,
\end{cases}
\end{equation}
where $\text{orth}(\cdot)$ denotes Gram-Schmidt orthonormalization.

Similar to HiSD, we incorporate the dimer method \eqref{the approximation of Hessian} to approximate $\mathbb{J}(\boldsymbol{x}) \boldsymbol{v}_i^{(n)}$ as
$$
\mathbb{J}(\boldsymbol{x}) \boldsymbol{v}_i^{(n)} \approx \dfrac{\boldsymbol{F}(\boldsymbol{x} + l \boldsymbol{v}_i^{(n)}) - \boldsymbol{F}(\boldsymbol{x} - l \boldsymbol{v}_i^{(n)})}{2l}.
$$
Combining this with the discretized $\boldsymbol{x}$-dynamics yields the complete discretization method:
\begin{equation}\label{discrete GHiSD}
\begin{cases}
\boldsymbol{x}^{(n+1)} = \boldsymbol{x}^{(n)} + \beta \left( \boldsymbol{F}(\boldsymbol{x}^{(n)}) - 2 \sum_{j=1}^{k} \left\langle \boldsymbol{F}(\boldsymbol{x}^{(n)}), \boldsymbol{v}_j^{(n)} \right\rangle \boldsymbol{v}_j^{(n)} \right), \\
\tilde{\boldsymbol{v}}_i^{(n+1)} = \boldsymbol{v}_i^{(n)} + \gamma \dfrac{\boldsymbol{F}(\boldsymbol{x}^{(n+1)} + l \boldsymbol{v}_i^{(n)}) - \boldsymbol{F}(\boldsymbol{x}^{(n+1)} - l \boldsymbol{v}_i^{(n)})}{2l}, \quad i = 1, \ldots, k, \\
\left[ \boldsymbol{v}_1^{(n+1)}, \ldots, \boldsymbol{v}_k^{(n+1)} \right] = \text{orth} \left( \tilde{\boldsymbol{v}}_1^{(n+1)}, \ldots, \tilde{\boldsymbol{v}}_k^{(n+1)} \right).
\end{cases}
\end{equation}

\subsection{AHiSD}\label{AHiSD: Accelerated High-index Saddle Dynamics}
This subsection introduces two momentum-accelerated variants of HiSD. The first variant employs heavy-ball momentum\cite{luo2025}:
\begin{equation}
\begin{cases}
\boldsymbol{x}^{(n+1)} = \boldsymbol{x}^{(n)} + \beta^{(n)} \left( \mathbb{I} - 2 \sum_{i=1}^{k} \boldsymbol{v}_i^{(n)} {\boldsymbol{v}_i^{(n)}}^\top \right) \boldsymbol{F}(\boldsymbol{x}^{(n)}) + \alpha (\boldsymbol{x}^{(n)} - \boldsymbol{x}^{(n-1)}), \\
\left\{\boldsymbol{v}_i^{(n+1)}\right\}_{i=1}^{k} = \text{EigenSol} \left( \mathbb{H}(\boldsymbol{x}^{(n+1)}), \left\{\boldsymbol{v}_i^{(n)}\right\}_{i=1}^{k} \right).
\end{cases}
\end{equation}

By utilizing the extrapolated iteration direction from momentum acceleration, we obtain the Nesterov acceleration variant \cite{Liu2024Neural}:
\begin{equation}
\begin{cases}
\boldsymbol{w}^{(n)} = \boldsymbol{x}^{(n)} + \alpha^{(n)} (\boldsymbol{x}^{(n)} - \boldsymbol{x}^{(n-1)}),\\
\boldsymbol{x}^{(n+1)} = \boldsymbol{w}^{(n)} + \beta^{(n)} \left( \mathbb{I} - 2 \sum_{i=1}^{k} \boldsymbol{v}_i^{(n)} {\boldsymbol{v}_i^{(n)}}^\top \right) \boldsymbol{F}(\boldsymbol{w}^{(n)}),\\
\left\{\boldsymbol{v}_i^{(n+1)}\right\}_{i=1}^{k} = \text{EigenSol} \left( \mathbb{H}(\boldsymbol{x}^{(n+1)}), \left\{\boldsymbol{v}_i^{(n)}\right\}_{i=1}^{k} \right).
\end{cases}
\end{equation}

Key features of these methods include:
\begin{itemize}
  \item \textbf{Eigensolver:} $\text{EigenSol}\big(\mathbb{H}(\boldsymbol{x}^{(n+1)}), \{\boldsymbol{v}_i^{(n)}\}_{i=1}^k\big)$ denotes the eigensolver, which may implement either direct HiSD dynamics discretization or algorithms such as LOBPCG.
  
  \item \textbf{Step Size:} The step size $\beta^{(n)}$ adopts the same selection criteria as HiSD, including Euler or BB step size methods.
  
  \item \textbf{Momentum:} The momentum term $\alpha (\boldsymbol{x}^{(n)} - \boldsymbol{x}^{(n-1)})$ constitutes the core acceleration mechanism. For heavy-ball acceleration, $\alpha \in [0,1)$ is typically constant. In Nesterov acceleration, $\alpha^{(n)}$ admits several formulations, such as
  \begin{equation}
      \alpha^{(n)} = \frac{n}{n+3},
  \end{equation}
  or
  \begin{equation}
      \alpha^{(n)} = \theta_{n+1}^{-1} (\theta_n - 1) \quad \text{with} \quad \theta_{n+1} = \frac{1 + \sqrt{1 + 4\theta_n^2}}{2},\; \theta_0 = 1.
  \end{equation}
  Numerical experiments show negligible performance differences between these choices; thus, we prefer the simpler first formulation.
\end{itemize}

\section{Construction of the Solution Landscape}

The ability to accurately locate individual index-$k$ saddle points using the HiSD algorithm suite mentioned above provides the foundation for a more ambitious goal: namely, mapping the global topology of the solution space. This section details the algorithmic framework for transforming discrete saddles into a pathway map consisting of all the stationary points and their interconnections, known as the \textbf{solution landscape}.

Conceptually, we model the solution landscape as a directed graph, G = (V, E), where:
\begin{itemize}
    \item Vertices (V) are the saddle points of the system, each characterized by its coordinates, Morse index and additional properties (e.g., their unstable manifolds).
    \item Edges (E) represent the hierarchical connections between these vertices. An edge from an index-$k$ saddle to an index-$m$ saddle (where $m<k$) signifies the latter saddle point can be reached from the unstable manifold of the former.
\end{itemize}

The construction of the solution landscape involves two primary operations, downward and upward searches, unified under the GHiSD framework, where gradient systems are subsumed by $\boldsymbol{F}(\boldsymbol{x}) = -\nabla E(\boldsymbol{x})$.

\subsection{Downward Search: Traversing the Hierarchy}
The downward search is the fundamental mechanism for unveiling the solution landscape. This search locates index-$m$ saddles ($m<k$) starting from index-$k$ saddles, terminating at index-$0$ saddles (sinks) while recording hierarchical connections. This process directly addresses a core challenge in landscape exploration: how to reliably escape a saddle point's local region and transition to connected, lower-index critical points. The methodology is as follows:

\begin{enumerate}
    \item \textbf{Initial perturbation:} Starting from a given index-$k$ saddle, a small perturbation of $\mathcal{O}(\varepsilon)$ is applied. To ensure an unbiased escape direction, a small random perturbation is first generated in the full space. This vector is then projected onto the unstable subspace $\mathcal{W}^{u}$ and scaled to $\mathcal{O}(\varepsilon)$ to extract the component that drives the system away from the saddle.
    \item \textbf{Unstable subspace initialization:} Select $m$ directions from the $k$ unstable directions at the initial saddle to form $\mathcal{V}_m^{(0)}$.
    \item \textbf{HiSD Algorithm:} The subsequent HiSD dynamics are applied, governed by $\boldsymbol{\dot{x}} = \left( \mathbb{I} - 2\mathcal{P}_{\mathcal{V}_m(\boldsymbol{x})} \right) \boldsymbol{F}(\boldsymbol{x})$, ensuring that the projection of the force along the unstable direction $\boldsymbol{v}_i$ remains positive. This drives the system away from the initial saddle and towards a connected, lower-index saddle, according to the theory of the HiSD algorithm. 
\end{enumerate}

To ensure a comprehensive exploration, the implementation incorporates enhancements such as sampling perturbations in the entire space from various statistical distributions (e.g., uniform distribution, Gaussian distribution) and systematically exploring all combinatorial possibilities ($k$ choose $m$) for the initial subspace. For sink searches ($m=0$), the algorithm elegantly simplifies. Since the subspace $\mathcal{V}^{(0)}$ is empty, the HiSD dynamics naturally reduce to the standard gradient descent dynamics, thereby bypassing expensive Hessian and eigenvector computations.

\subsection{Upward Search: Exploring from the Basins}
While downward search builds the primary hierarchical structure, an upward search provides a complementary strategy for exploring the solution landscape, enabling discovery of higher-index saddles from lower-index ones. This is particularly useful for initiating new searches from discovered minima or expanding an incomplete landscape. The principles are symmetric to the downward search:

\begin{enumerate}
    \item \textbf{Initial perturbation:} From the initial index-$k$ saddle, apply a small perturbation of $\mathcal{O}(\varepsilon)$. This perturbation is specified by user when restarting.
    \item \textbf{Subspace initialization:} Form $\mathcal{V}_m^{(0)}$ by combining the $k$ unstable directions with $m-k$ directions from the stable subspace (typically those corresponding to the smallest positive eigenvalues are chosen).
    \item \textbf{HiSD Algorithm:} The subsequent HiSD dynamics induce a negative projection along the stable direction $\boldsymbol{v}_{j}$, enabling escape.
\end{enumerate}

In our practical implementation, the overall strategy begins by using HiSD to locate an initial high-index saddle, which then initiates a comprehensive downward search to construct the solution landscape. Upward searches are also valuable tools for restarting the exploration from newly discovered regions to ensure completeness.

\subsection{Algorithm Framework for Solution Landscape Construction}

When constructing the entire solution landscape, there are some key challenges beyond escaping problems.
\begin{itemize}
    \item \textbf{The identification problem:} The search process will frequently converge to the same saddle point from different paths. A robust and efficient method is needed to determine if a newly found critical point is genuinely novel or a duplicate of an existing one in the list.
    \item \textbf{The completeness problem:} The search must be systematic and exhaustive enough to discover almost all significant connections and build a nearly complete representation of the solution landscape, avoiding premature termination down a single path.
\end{itemize}

To address these challenges, we designed a global search strategy based on the BFS algorithm, managed with a FIFO queue. We chose BFS over alternatives like Depth-First Search (DFS) for a critical reason: its level-by-level exploration perfectly mirrors the natural hierarchy of the solution landscape. Whereas DFS would rapidly follow one path to a minimum, potentially missing other important branches, BFS guarantees that all index-$(k-1)$ saddles connected to a specific index-$k$ saddle are found before the search proceeds to the index-$(k-2)$ level from the same initial saddle. This methodical approach is essential for ensuring both systematicity and completeness.

The complete automated workflow is detailed in the following algorithm.

\begin{algorithm}[H]
\caption{\textbf{Solution Landscape Construction via HiSD Search.}}
\label{alg:landscape}

\textbf{Input:} Initial point $\boldsymbol{x^{(0)}}$, initial unstable subspace $\boldsymbol{E^{(0)}}$(can be constructed automatically), maximum Morse index $M$, search depth for one point $J$

\begin{algorithmic}[1]
\State Perturb initial point to avoid being trapped in the initial saddle point: $\boldsymbol{x^{(0)}} \gets \boldsymbol{x^{(0)}} + \boldsymbol{x_{\epsilon}}$
\State Initialize the HiSD solver
\For{$i=M$ downto $0$}
    \State Launch HiSD search for index-$i$ saddle points using $\boldsymbol{x^{(0)}}$ and $\boldsymbol{E^{(0)}}$
    \If{New saddle point $\boldsymbol{x_{\textbf{new}}}$ found}
        \State Add $\boldsymbol{x_{\textbf{new}}}$ to the saddle-point list and initialize its connection data
        \State \textbf{break}
    \EndIf
\EndFor  
\State \textit{// Implement BFS with FIFO queue scheduling}
\While{Queue non-empty}
    \State Dequeue head element $\boldsymbol{x_{\textbf{head}}}$ (index-$k$)
    \For{$i=k-1$ downto $\max(0,k-J)$}
        \State Generate symmetric perturbation set $\boldsymbol{P}$ (uniform/Gaussian)
        \State Generate initial subspaces $\boldsymbol{E}$ when $i>0$ (No need for $i=0$)
            \Statex \qquad \quad \textit{// Try all combinations $\binom{k}{i}$ of unstable directions if parameter is 'all',}
\Statex \qquad \quad \textit{// otherwise select the $i$ most unstable directions if parameter is 'min'.}
            \ForAll{$\boldsymbol{p} \in \boldsymbol{P}, \boldsymbol{E_i} \in \boldsymbol{E}$}
                \State Launch HiSD search from $\boldsymbol{x_{\textbf{head}}}+\boldsymbol{p}$ with $\boldsymbol{E_i}$
                \If{New saddle $\boldsymbol{x_{\textbf{new}}}$ found}
                    \State Compare with existing points
                    \If{$\boldsymbol{x_{\textbf{new}}}$ is novel}
                        \State Enqueue and update connection info
                    \Else
                        \State Record connection with $\boldsymbol{x_{\textbf{head}}}$
                    \EndIf
                \EndIf
            \EndFor
    \EndFor
\EndWhile    
\State (Optional) Repeat process from new initial point
\end{algorithmic}
\textbf{Output:} The saddle-point list (IDs, coordinates, Morse indices, interconnections)
\end{algorithm}

\section{Software Implementation}

\subsection{Overview}
This section systematically presents the functionality of each component within the software package and demonstrates implementation details of the core code. We begin by outlining Python package dependencies and input parameter specifications.

\subsubsection{Dependencies}
Our package requires the following Python libraries under stable versions:
\begin{itemize}
    \item \textbf{Standard Library Modules:} \graycode{copy}, \graycode{sys}, \graycode{warnings}, \graycode{inspect}, \graycode{json}, \graycode{itertools}, \graycode{math}, and \graycode{pickle}
    \item \textbf{External Dependencies:} \graycode{numpy} (v2.2.3), \graycode{scipy} (v1.15.3), \graycode{sympy} (v1.13.3), \graycode{matplotlib} (v3.10.0), and \graycode{networkx} (v3.4.2).
\end{itemize}

\subsubsection{Input Parameter Configuration}\label{parameter}
This section provides an overview of the available parameters that offer flexible and customizable solvers. For a detailed parameter list, refer to Appendix \ref{Parameter Reference}.
\begin{itemize}

\item \textbf{System Parameters:} These parameters are related to the system setup and its general properties, including \textbf{Dim}, \textbf{EnergyFunction}, \textbf{Grad}, \textbf{AutoDiff}, \textbf{NumericalGrad}, \textbf{DimerLength}, \textbf{SymmetryCheck}, and \textbf{GradientSystem}.

\textbf{Note:} To maintain formulaic uniformity between gradient systems (governed by $-\nabla E$) and non-gradient systems (governed by $\boldsymbol{F}$), the package employs an internal sign convention. Specifically, when providing input for non-gradient systems, users must supply the negative of the vector field $\boldsymbol{F}(\boldsymbol{x})$. This ensures that the core iterative machinery of HiSD and GHiSD remains consistent regardless of the system type.

\item \textbf{Landscape Parameters:} These parameters are related to constructing and navigating the solution landscape, including \textbf{MaxIndex}, \textbf{MaxIndexGap}, \textbf{SameJudgement}, \textbf{InitialEigenVectors}, \textbf{PerturbationMethod}, \textbf{PerturbationRadius}, \textbf{PerturbationNumber}, and \textbf{EigenCombination}. They govern the hierarchical search for saddle points, including index limits, perturbation strategies, and equivalence criteria.

\item \textbf{Solver Parameters:} These parameters are related to the solver process and control the behavior of the HiSD algorithm, including \textbf{InitialPoint}, \textbf{Tolerance}, \textbf{SearchArea}, \textbf{TimeStep}, \textbf{MaxIter}, \textbf{SaveTrajectory}, \textbf{Verbose}, and \textbf{ReportInterval}. They define the initial conditions, convergence criteria, and output settings for saddle point search.

\item \textbf{Eigen Parameters:} These parameters control how the eigenpairs (eigenvalues and eigenvectors) are computed, including \textbf{EigenMethod}, \textbf{EigenMaxIter}, \textbf{EigenStepSize}, \textbf{PrecisionTol}, and \textbf{EigvecUnified}. They define the eigensolver selection, iteration limits, step sizes, and tolerance for eigenvalue classification.

\item \textbf{Hessian Parameters:} These parameters are related to the Hessian matrix, including \textbf{ExactHessian} and \textbf{HessianDimerLength}. They control the computation method for Hessian matrices, such as analytical derivation or numerical approximation via dimer-based approaches.

\item \textbf{Acceleration Parameters:} These parameters are related to improving the speed and efficiency of the algorithm, including \textbf{BBStep}, \textbf{Acceleration}, \textbf{NesterovChoice}, \textbf{NesterovRestart}, and \textbf{Momentum}. They enable adaptive step sizing and momentum-based acceleration techniques to enhance convergence.

\end{itemize}

\subsection{Core Program \graycode{SaddleScape.py}}
This program encapsulates all functionalities of the software package within the \graycode{Landscape} class, which primarily contains three types of functions: 
\begin{itemize}
\item \textbf{\graycode{\_\_init\_\_}:} The essential function for class definition, performing parameter validation and initialization (see Section \ref{Parameter Validation Program}).
\item Functions for constructing the solution landscape using saddle point dynamics (see Section \ref{Solution Landscape Run Program}):
    \begin{itemize}
      \item \textbf{\graycode{Run}:} The basic solver.
      \item \textbf{\graycode{RestartFromPoint}:} A function that restarts the computation from a new initial point.
      \item \textbf{\graycode{RestartFromSaddle}:} A function that restarts the computation from an existing saddle point in the solution landscape, using a small random perturbation to avoid being trapped near the original saddle point.
    \end{itemize}
\item Functions for visualization and saving:
\begin{itemize}
  \item \textbf{\graycode{DrawTrajectory}:} For plotting solution trajectories (see Section \ref{Search Trajectory Plotting Program}).
  \item \textbf{\graycode{DrawConnection}:} For visualizing the solution landscape by drawing saddle point interconnections (see Section \ref{Solution Landscape Plotting Program}).
  \item \textbf{\graycode{Save}:} For storing raw computational data (see Section \ref{Data Saving Program}).
\end{itemize}
\end{itemize}

\subsection{Parameter Validation Module}\label{Parameter Validation Program}

\subsubsection{Parameter Validation Program \graycode{LandscapeCheckParam.py}}

This program primarily validates the \textbf{Landscape Parameters} in user inputs to ensure compliance with requirements via the function \graycode{LandscapeCheckParam}. Its core functions include:

\begin{itemize}

\item Raising immediate errors when invalid inputs are detected.

\item Automatically assigning default values for missing parameters.

\end{itemize}

\subsubsection{Parameter Validation Program \graycode{HiSDSolverCheckParam.py}}

The remaining parameters are validated by calling the function \graycode{HiSDCheckParam} from the file \graycode{HiSDSolverCheckParam.py} when \graycode{instance.calHiSD} is instantiated using the \graycode{Solver} class. Its implementation logic and functionality are similar to \graycode{LandscapeCheckParam.py}.

\subsection{Architecture and Implementation of Index-$k$ Saddle Point Search \graycode{HiSDSolver.py}}\label{HiSD Run Program}
The solver architecture is designed to accommodate the diverse configurations detailed in Appendix \ref{Parameter Reference}. This section details how the index-$k$ saddle point search program addresses these diverse needs, including:  
\begin{itemize}
  \item Comprehensive descriptions of algorithmic branches under different parameter configurations, including functional specifications and recommended usage scenarios;
  \item Concrete implementation details of the core saddle point search algorithm.
\end{itemize}

The final integrated solver with specified configurations is encapsulated within the \graycode{Solver} class:
\begin{itemize}
  \item \textbf{\graycode{\_\_init\_\_}:} Performs parameter validation and initialization (via the \graycode{HiSDCheckParam} function described above).
  \item \textbf{\graycode{run}:} Executes index-$k$ saddle point search when called (the target search index is specified during \graycode{Solver} initialization).
\end{itemize}

\subsubsection{Symbolic Computation and Automatic Differentiation \graycode{AnaCal.py}}
This program implements a dual interface design for energy/gradient inputs with automated symbolic-numeric conversion:

\begin{enumerate}
    \item \textbf{Direct use of callable objects:} When function objects are provided as energy functions and gradient functions, they are used directly without additional processing (No specific function call required for this case).
    
    \item \textbf{Symbolic energy function processing:}
    \begin{itemize}
        \item \textbf{\graycode{EnergyFunctionAnalysis}:} Converts symbolic energy functions into callable numerical functions, avoiding repeated compilation of symbolic expressions to significantly improve computational efficiency. Results are stored in \graycode{instance.EnergyFunctionSolve}.
        \item \textbf{\graycode{EnergyFunctionCalculate}:} Evaluates energy functions at specific points using the function compiled by \graycode{EnergyFunctionAnalysis}.
    \end{itemize}
    
    \item \textbf{Symbolic gradient processing:}
    \begin{itemize}
        \item \textbf{\graycode{AutoDerivative}:} Performs automatic differentiation of symbolic energy functions, returning symbolic differentiation results and callable numerical gradient functions. The symbolic results are stored in \graycode{instance.GradFunction\_sympy}, while the numerical functions reside in \graycode{instance.GradFunction\_numpy}.
        \item \textbf{\graycode{ExactGradAnalysis}:} Compiles user-provided symbolic gradient functions into callable numerical functions, with output and storage methods identical to \graycode{AutoDerivative}.
        \item \textbf{\graycode{AutoGrad}:} Computes gradients at specific points using pre-stored compiled gradient functions in \graycode{instance.GradFunction\_numpy}.
    \end{itemize}
    
    \item \textbf{Numerical gradient computation:}
    \begin{itemize}
        \item \textbf{\graycode{AutoGradNum}:} Handles cases where only numerical energy functions are provided, computing numerical approximations of gradient functions via central difference quotients.
    \end{itemize}
\end{enumerate}

\subsubsection{Hessian Matrix Processing \graycode{HessianMatrix.py}}
\noindent\textbf{Note:} For non-gradient systems, the term "Hessian" here refers to the \textbf{negative} Jacobian matrix of the vector field function. This is because we previously required users to prepend a negative sign to the force input.

\begin{enumerate}
    \item \textbf{Treating Hessian as an operator:} In the HiSD algorithm, we only need to handle Hessian-vector or Hessian-matrix structures, avoiding explicit computation of the full Hessian matrix.
    \begin{itemize}
        \item \textbf{\graycode{SingleHessianVectorProduct}:} Computes Hessian-vector products using central finite difference approximation with dimers $\mathbf{x}\pm L\mathbf{v}$ (implements Eq. \eqref{the approximation of Hessian}).
        \item \textbf{\graycode{BatchHessianVectorProduct}:} Computes Hessian-matrix products using matrix-free methods via Hessian-vector operator overloading.
    \end{itemize}
    
    \item \textbf{Explicit Hessian computation:}
    \begin{itemize}
        \item \textbf{\graycode{Hessian\_Analysis\_withIfsym}:} Performs symbolic differentiation with symmetry checks to determine if the system is gradient-preserving, which determines the subsequent HiSD algorithms to be used.
        \item \textbf{\graycode{Hessian\_Analysis}:} Computes Hessian matrices without symmetry checks.
        \item \textbf{\graycode{Hessian\_exact}:} Evaluates Hessian matrices at specific points using symbolic results.
        \item \textbf{\graycode{Hessian}:} Hybrid approach for explicit Hessian computation:
        \begin{itemize}
            \item Uses exact symbolic computation when specified.
            \item Otherwise employs numerical approximation via \graycode{BatchHessianVectorProduct} (Hessian times columns of the identity matrix, i.e., reconstructing the full Hessian by applying the operator to each standard basis vector).
        \end{itemize}
        \item \textbf{Applicable scenarios for explicit Hessian computation:}
        \begin{itemize}
            \item When iteration counts may be very large.
            \item When Hessian reuse is frequent.
            \item When solution accuracy takes priority over computation time (in which case only symbolic explicit computation is used).
        \end{itemize}
    \end{itemize}
\end{enumerate}

\textbf{Note}: Explicit numerical computation of Hessian remains discouraged in high-dimensional cases due to its substantial computational cost.

\subsubsection{Eigenvalue Solver \graycode{EigMethod.py}}
This program implements various methods to calculate eigenpairs for HiSD algorithms, supporting both gradient and non-gradient systems.

\begin{enumerate}
    \item \textbf{LOBPCG Method:} (exclusively for gradient systems due to its requirement for symmetric operators) This method is implemented using \graycode{scipy.sparse.linalg.lobpcg}, corresponding to Section \ref{Discretization of HiSD method}.
    \begin{itemize}
        \item \textbf{\graycode{lobpcg}:} Implements LOBPCG using matrix-free operators (Hessian-vector and Hessian-matrix products).
        \item \textbf{\graycode{lobpcg\_exacthessian}:} Uses explicit Hessian matrices for LOBPCG computation when exact Hessian is available.
        \item Both functions return:
        \begin{itemize}
            \item Eigenvectors corresponding to the smallest $k$ eigenvalues at the current point.
            \item A flag \graycode{whetherkindex} indicating whether the point satisfies the index-$k$ condition (at least $k$ negative eigenvalues). Since output eigenvalues are sorted in ascending order, we only need to check if the $k$-th eigenvalue is negative (under the precision tolerance).
        \end{itemize}
    \end{itemize}

    \item \textbf{Power Method:} (for non-gradient systems) This method is based on direct discretization of the theoretical dynamical system $\boldsymbol{\dot{v}} = \mathbb{J}(\boldsymbol{x}) \boldsymbol{v}$ (see Eq. \eqref{discrete GHiSD of W^u} in Section \ref{Dis_GHiSD}). After each update of the approximate eigenvectors, QR decomposition is applied for orthonormalization.
    \begin{itemize}
        \item \textbf{\graycode{power\_nonGrad}:} Matrix-free implementation using Hessian-vector and Hessian-matrix products.
        \item \textbf{\graycode{power\_nonGrad\_exacthessian}:} Explicit matrix implementation when exact Hessian is available.
        \item Outputs are consistent with LOBPCG. To reduce computational cost from repeated exact eigenvalue solutions, we approximate using diagonal elements of $\mathbb{V}^{\top}\mathbb{H}\mathbb{V}$, where:
            \begin{itemize}
                \item $\mathbb{H}$ is the Hessian (here,the negative Jacobian) matrix at the point.
                \item $\mathbb{V}$ contains orthogonal approximate unit eigenvectors as its columns.
            \end{itemize}
    \end{itemize}
    
    \textbf{Note:} In fact, since gradient systems are essentially a special case of non-gradient systems, the power method can also be applied to gradient systems, provided that the interface remains consistent.

    \item \textbf{Explicit Euler Discretization:} (for both gradient and non-gradient systems) This method implements direct Euler discretization of the continuous HiSD and GHiSD dynamics from the theoretical framework. For further details, please refer to Algorithm \ref{alg:hisd} in Section \ref{Discretization of HiSD method} and Eq. \eqref{discrete GHiSD} in Section \ref{Dis_GHiSD}.
    \begin{itemize}
        \item \textbf{\graycode{euler}:} Matrix-free Euler discretization for continuous HiSD dynamics in gradient systems.
        \item \textbf{\graycode{euler\_exacthessian}:} Explicit Hessian matrix discretization for gradient systems.
        \item \textbf{\graycode{euler\_nonGrad}:} Matrix-free Euler discretization for continuous GHiSD dynamics.
        \item \textbf{\graycode{euler\_nonGrad\_exacthessian}:} Euler discretization for continuous GHiSD dynamics.
        \item Output and eigenvalue approximation approaches remain consistent with the Power Method.
    \end{itemize}

    \textbf{Note:} Considering the equivalence between direct explicit Euler discretization and the power method as discussed in Section \ref{Dis_GHiSD}, and for computational efficiency, the software package adopts the power method as the default eigenvalue method for non-gradient systems, and LOBPCG, which is specially optimized for gradient systems, as the default for gradient systems. However, if users are more concerned about the theoretical guarantees of HiSD, they can switch to using direct explicit Euler discretization.
    
    \item \textbf{Saddle Point Index Verification Programs:} Provide accurate saddle point index verification.
    \begin{itemize}
        \item \textbf{\graycode{CheckIndexk}:} Provides a precise sparse method to calculate the smallest $k$ eigenvalues, avoiding inaccuracies from approximation methods. This function is specifically called when HiSD iteration approaches a new saddle point instead of every iteration to reduce computational cost.
        \item \textbf{\graycode{FindIndex}:} Rechecks the Morse index after iteration completes (since termination conditions only guarantee at least index-$k$ due to accuracy limitations). It computes the full spectrum of the Hessian and treats eigenvalues with sufficiently small real parts as approximate zeros (accounting for numerical errors), helping detect potential saddle point degeneracy.
    \end{itemize}
    
    \item \textbf{Initial Eigenvector Selection Program:} Provides approximate eigenvectors for initial iteration points when no initial eigenvectors are given.
    \begin{itemize}
        \item \textbf{\graycode{GiveEigenvector}:} Applies symmetrization uniformly to both gradient and non-gradient systems to accelerate computation, as the algorithm's minimum requirement is simply obtaining a set of orthogonal unit vectors to serve as initial eigenvectors for iteration.
    \end{itemize}
    
    \item \textbf{Canonicalization Program for Eigenvectors:} Due to significant differences in the implementation of fundamental algorithms (e.g., LAPACK) across hardware platforms, eigenvectors for multiple eigenvalues may vary substantially between devices, even though they all form valid bases for the same eigenspace. When set to \graycode{True}, this parameter activates a canonicalization procedure (Gaussian elimination to Row Echelon Form followed by Gram-Schmidt orthogonalization) to standardize the basis vectors of the eigenspace, thereby improving cross-device reproducibility.

    In detail, the same standardization algorithm is applied to each eigenblock based on eigenvalue multiplicity. This procedure is only needed for complicated systems exhibiting multiple eigenvalues. In cases of simple multiplicity, only negligible signature discrepancies may arise between devices. As demonstrated by Eq. \eqref{2.4} and Eq. \eqref{Entire GHiSD}, such minor variations do not affect the convergence behavior of the solution iterations.

    \textbf{Note:} Outputs may still vary between devices due to differences in rounding-error propagation.
    \begin{itemize}
        \item \textbf{\graycode{CanonicalizeEigens}:} Standardizes eigenvectors for multiple eigenvalues to improve cross-device stability.
        \item \textbf{\graycode{CanonicalizeMatrix}:} Canonicalizes row space using Gaussian elimination and Gram-Schmidt orthogonalization.
    \end{itemize}
\end{enumerate}

\subsubsection{Iterative Search Program \graycode{HiSDIter.py}}

Now, we present the key part of the solver: the search iterations.

\noindent\textbf{Preparation Programs}

\begin{itemize}
    \item \textbf{\graycode{HiSDInitialization}:} Initializes the initial point and initial eigenvectors required for iteration.
    \item \textbf{\graycode{BBStepSize}:} Implements the Barzilai-Borwein adaptive step size strategy based on Eq. \eqref{BB2step}.
\end{itemize}

\noindent\textbf{\graycode{HiSDIteration}: For Searching Non-index-$0$ Saddle Points}

With the preparatory steps completed, the program structure is as follows. This section mainly focuses on:
\begin{enumerate}
    \item \textbf{Momentum-based acceleration schemes:}
        \begin{itemize}
            \item Heavy-ball acceleration (including the special case of no acceleration, equivalent to heavy-ball with momentum coefficient $0.0$).
            \item Nesterov acceleration algorithm.
        \end{itemize}
    \item \textbf{Step size selection strategies:}
        \begin{itemize}
            \item Use fixed time step \graycode{dt=instance.TimeStep}.
            \item Adopt BB step size \graycode{dt=BBStepSize(x, x\_pre, g, g\_pre, tau)}.
        \end{itemize}
    \item \textbf{Constructing the complete iteration framework:} Calling appropriate functions according to the parameter configuration.
    \item \textbf{Returning newly found saddle point information and iteration records:} Logging the iteration count and gradient norm every \graycode{instance.ReportInterval} steps when \graycode{instance.Verbose=True}.
    \item \textbf{Warning and error handling during iterations:}
    \begin{itemize}
        \item Issues warning when search exceeds defined boundaries or fails to converge within maximum iterations, without updating results in \graycode{instance}.
        \item Issues warning for potentially degenerate saddle points (while still recording results in \graycode{instance}), reports the index of found saddle point for non-degenerate cases.
    \end{itemize}
\end{enumerate}
\begin{lstlisting}[caption={\textbf{Implementation of the core iterative solver for HiSD,  including selectable momentum acceleration.}}, label=HiSDIter.py(2), language=Python]
def HiSDIteration(instance):
	"""
	Main iteration for HiSD method for high-index saddle points.
	"""
	boundarytempsave = copy.deepcopy(instance.DataBoundary)
	xini = instance.PrimaryInitialPoint.reshape(
		-1,
	)
	tau = 0.5
	x = instance.x
	x_pre = instance.x_pre
	v = instance.v
	WhetherBB = instance.BBStep
	mom = instance.Momentum
	dt = instance.TimeStep
	if instance.SaveTrajectory:
		x_record = copy.deepcopy(x.reshape(1, instance.Dim))
	else:
		x_record = None
	gnorm_record = []
	timestep_record = [0]
	g = instance.Grad(x)
	gnorm_record.append(np.linalg.norm(g))

	# Check if Nesterov acceleration is used
	if instance.Acceleration == "nesterov":
		if instance.NesterovChoice == 2:
			NesterovTheta = 1
			NesterovTheta = (1 + (1 + 4 * NesterovTheta**2) ** (1 / 2)) / 2
		elif instance.NesterovChoice == 1:
			if instance.NesterovRestart is not None:
				plusbase = 0

		# Main iteration loop
		for j in range(1, instance.MaxIter + 1):
			if WhetherBB and j > 1:
				dt = BBStepSize(x, x_pre, g, g_pre, tau)
			timestep_record.append(dt)
			if instance.NesterovChoice == 2:
				NesterovThetaNew = (1 + (1 + 4 * NesterovTheta**2) ** (1 / 2)) / 2
				NesterovGamma = (NesterovTheta - 1) / NesterovThetaNew
				NesterovTheta = NesterovThetaNew
			elif instance.NesterovChoice == 1:
				if instance.NesterovRestart is not None:
					NesterovGamma = (j - plusbase) / (j - plusbase + 3)
				else:
					NesterovGamma = j / (j + 3)
			w_temp = x + NesterovGamma * (x - x_pre)
			gw = instance.Grad(w_temp)
			dxw = dt * (gw - 2.0 * np.matmul(v, np.matmul(v.T, gw)))
			x_temp = w_temp - dxw
			x_pre = x
			x = x_temp
			if WhetherBB:
				g_pre = g
			v, whetherkindex = instance.EigVecMethod(x, v)
			g = instance.Grad(x)
			if instance.Verbose:
				if j % instance.ReportInterval == 0:
					print("Iteration: "+ str(j)+ f"|| Norm of gradient: {np.linalg.norm(g):.6f}")
			xnow = x.reshape(
				-1,
			)
			if instance.SaveTrajectory:
				x_record = np.append(x_record, x.reshape(1, instance.Dim), axis=0)
			if np.linalg.norm(xnow - xini) > instance.SearchArea:
				print("[WARNING] Iteration diverged: Search point exceeds feasible region. Skipping to next search.")
				instance.flag = False
				instance.DataBoundary = boundarytempsave
				return
			instance.DataBoundary = [[min(instance.DataBoundary[i][0], x[i, 0]),max(instance.DataBoundary[i][1], x[i, 0]),]for i in range(instance.Dim)]
			gnorm_record.append(np.linalg.norm(g))
			if np.linalg.norm(g) < instance.Tolerance:
				if whetherkindex:
					break
				if j==1 or gnorm_record[-2]>=instance.Tolerance:
					whetherk = CheckIndexk(instance, x, instance.SaddleIndex)
					if whetherk:
						break
			if instance.NesterovRestart is not None:
				if j % instance.NesterovRestart == 0:
					if instance.NesterovChoice == 2:
						NesterovTheta = 1
						NesterovTheta = (1 + (1 + 4 * NesterovTheta**2) ** (1 / 2)) / 2
					elif instance.NesterovChoice == 1:
						plusbase = j
	else:
		... # Standard/Heavy-ball iteration logic goes here

	# Check for convergence and report results
	if j == instance.MaxIter:
		if gnorm_record[-1] >= instance.Tolerance:
			print("[WARNING] Iteration not converged: Maximum iterations reached without convergence. Skipping to next search.")
			instance.flag = False
			instance.DataBoundary = boundarytempsave
			return
		negnum, zeronum, posnum, negvecs = FindIndex(instance, x)
		if negnum+zeronum >= instance.SaddleIndex:
			print("[Note] Due to eigenvalue approximation inaccuracies during iterations, "
				"the search trajectory may reach a qualifying saddle point without triggering a report.")
		else:
			print("[WARNING] Iteration not converged: Maximum iterations reached without convergence. Skipping to next search.")
			instance.flag = False
			instance.DataBoundary = boundarytempsave
			return
	else:
		negnum, zeronum, posnum, negvecs = FindIndex(instance, x)
	if zeronum != 0:
		print(f"[WARNING] Degenerate saddle point detected under precision tol={instance.PrecisionTol}: Hessian matrix may contain zero eigenvalue(s).")
		print(f"Eigenvalue spectrum: negative={negnum}, zero={zeronum}, positive={posnum}. ")
	else:
		print(f"Non-degenerate saddle point identified: Morse index ={negnum} (number of negative eigenvalues).")
	instance.finalindex = int(negnum)
	instance.negvecs = negvecs
	instance.IterNum = j
	instance.x = x
	instance.x_record = x_record
	instance.timestep_record = np.array(timestep_record).reshape(-1, 1)
	instance.timestep_record = np.cumsum(instance.timestep_record)
	instance.gnorm_record = gnorm_record
\end{lstlisting}

\noindent\textbf{\graycode{SDIteration}: For Searching Index-$0$ Saddle Points}

This function is a simplified version of \graycode{HiSDIteration} specifically for index-$0$ saddle points (minima), noting that no eigenvectors are used during iterations. It is based on gradient descent (optionally with acceleration algorithms) and does not require Hessian matrix information to separate eigen-directions (e.g., no need to compute eigenvectors), thereby significantly reducing unnecessary computations.

\subsubsection{Summary}
The implementation consists of preparatory modules for: (1) symbolic computation and automatic differentiation, (2) Hessian matrix processing, and (3) eigenvalue solving. These are integrated into the parameter validation program (\graycode{HiSDSolverCheckParam.py}) and the iterative search program (\graycode{HiSDIter.py}). The final integrated \graycode{HiSDSolver.py} serves as the core module for subsequent solution landscape construction.

\subsection{Solution Landscape Construction \graycode{LandscapeRun.py}}\label{Solution Landscape Run Program}

This program serves as the central orchestration module, responsible for constructing the solution landscape. All previous components are designed to achieve this objective. Given the program's length, we will adopt a layered decomposition methodology.

\subsubsection{Preparation: Function \graycode{checkwhetherexist}}
\begin{lstlisting}[caption={\textbf{Function for uniquely identifying saddle points using user-defined equivalence metrics to avoid duplicates.}}, label=LandscapeRun.py(1), language=Python]
def checkwhetherexist(instance, x, index):
	"""
	Check if a saddle point already exists.
	"""
	for i in range(len(instance.SaddleList)):
		if instance.SameJudgementMethod(x, instance.SaddleList[i][1]) and (index == instance.SaddleList[i][2]):
			return True, i
	return False, None
\end{lstlisting}

This function is used to distinguish distinct saddle points and merge duplicates during upward and downward searches.

\subsubsection{Core Code: Function \graycode{LandscapeRun}}

\noindent\textbf{Searching for an Index-\graycode{MaxIndex} Saddle Point from the Initial Point}

This part of the algorithm aims to find a saddle point, starting from the given initial point, having the highest possible index not exceeding \graycode{instance.MaxIndex}.

The found saddle point serves as the starting point for the downward search. The loop terminates immediately upon finding a suitable saddle point; otherwise, the procedure reports failure.

\textbf{Note:} As this function is also used in restart programs \graycode{RestartFromPoint} and \graycode{RestartFromSaddle}, additional requirements include:
\begin{itemize}
  \item The point should satisfy the separation condition from existing saddle points (checked by \graycode{LandscapeRun}).
  \item For points deemed identical, the one with smaller gradient norm is kept.
\end{itemize}
These proximity checks are skipped during initial landscape construction when no saddle-point list exists.
\begin{lstlisting}[caption={\textbf{Algorithm to locate the initial saddle point when beginning or restarting.} The routine attempts to find a saddle of the specified \texttt{MaxIndex} from the initial guess, iteratively decrementing the target index if convergence fails, to seed the solution landscape.}, label=LandscapeRun.py(2), language=Python]
def LandscapeRun(instance):
	"""
	Do the downward search automatically of HiSD to construct the solution landscape
	"""
	tempList = []  # Temporary list to hold intermediate saddle points
	k = instance.MaxIndex
	np.random.seed(1121)
	# Searching for saddle points starting from MaxIndex
	while k >= 0:
		instance.calHiSD.SaddleIndex = k
		if instance.InitialEigenVectors is not None:
			instance.calHiSD.InitialEigenVectors = instance.InitialEigenVectors[:, 0:k]
		... # Formatted output
		instance.calHiSD.run()
		if instance.calHiSD.flag:
			whetherinlist, smallestind = checkwhetherexist(instance, instance.calHiSD.x, instance.calHiSD.finalindex)
			if whetherinlist:
				xold = instance.SaddleList[smallestind][1]
				xnew = instance.calHiSD.x
				gold = instance.calHiSD.Grad(xold)
				gnew = instance.calHiSD.Grad(xnew)
				if np.linalg.norm(gold) > np.linalg.norm(gnew):
					instance.SaddleList[smallestind][1] = instance.calHiSD.x
				if instance.Continue:
					tempList.append([instance.SaddleList[smallestind][0], instance.SaddleList[smallestind][1], 
							instance.SaddleList[smallestind][2], instance.SaddleList[smallestind][3]])
				if instance.BeginID != -1 and (
					instance.BeginID not in instance.SaddleList[smallestind][4]
				):
					instance.SaddleList[smallestind][4].append(instance.BeginID)
					if instance.SaveTrajectory:
						instance.DetailRecord.append([smallestind, instance.BeginID, 
								instance.calHiSD.x_record, instance.calHiSD.timestep_record])
					else:
						instance.DetailRecord.append([smallestind, instance.BeginID])
				print("Search an existing saddle point.")
			else:  # Save the saddle point and its information
				instance.SaddleList.append([instance.saddleind, instance.calHiSD.x, 
						instance.calHiSD.finalindex, instance.calHiSD.negvecs, [instance.BeginID]])
				tempList.append([instance.saddleind, instance.calHiSD.x,
						instance.calHiSD.finalindex, instance.calHiSD.negvecs])
				if instance.SaveTrajectory:
					instance.DetailRecord.append([instance.saddleind, instance.BeginID,
							instance.calHiSD.x_record, instance.calHiSD.timestep_record])
				else:
					instance.DetailRecord.append([instance.saddleind, instance.BeginID])
				instance.saddleind = instance.saddleind + 1
			if instance.calHiSD.finalindex > instance.MaxIndex:
				instance.MaxIndex = instance.calHiSD.finalindex
				print("Warning! 'MaxIndex' updated due to a larger saddle point index.")
			break
		k = k - 1
	if k == -1:
		raise RuntimeError("No more saddle points found in the search area!")
\end{lstlisting}

\noindent\textbf{Breadth-First Downward Search from High-index Saddle Points}

The core strategy for locating index-$j$ ($j < k$) saddle points from an index-$k$ saddle point involves the following steps:

\begin{enumerate}
\item \textbf{Direction Selection:} Select $j$ unstable directions from the $k$ available modes to form the maintained subspace. These unstable directions correspond to negative eigendirections of the Hessian in gradient systems, or positive real-part eigendirections of the Jacobian for non-gradient systems; note that our convention unifies these by treating the latter as negative directions of the 'Hessian'.
    
    \item \textbf{Perturbation:} Apply perturbations to the coordinates of the current saddle point.
    
    \item \textbf{Escape Mechanism:} Once the perturbation displaces the state from the equilibrium along the remaining $k-j$ unstable directions:
    \begin{itemize}
        \item The system forces drive the dynamics along these $k-j$ directions.
        \item Consequently, $\boldsymbol{x}$ escapes the neighborhood of the current index-$k$ saddle.
        \item The HiSD iterations then guide the evolution toward a target index-$j$ saddle point.
    \end{itemize}
\end{enumerate}

To optimize computational efficiency, the implementation distinguishes between general saddle point searches ($j>0$) and local minima searches ($j=0$):

\begin{itemize}
    \item \textbf{For $j > 0$ cases:}
    \begin{itemize}
        \item Initialize the unstable subspace with a combination of $j$ eigenvectors selected from the unstable modes (using either all combinatorial subsets or the most unstable modes).
        \item Generate initial random perturbations (e.g., following a Uniform or Gaussian distribution).
        \item Proceed with the standard HiSD search and update the data structures upon convergence.
    \end{itemize}

    \item \textbf{For $j = 0$ cases:}
    \begin{itemize}
        \item The eigen-subspace analysis is omitted, and a direct gradient descent method is applied. This simplification maintains search effectiveness for local minima while significantly reducing computational costs.
    \end{itemize}
\end{itemize}

\textbf{Note:} Considering the approximate symmetry of the solution landscape, we apply perturbations in opposite directions (i.e., $\pm \boldsymbol{p}$) to identify connected saddle points more comprehensively.

\begin{lstlisting}[caption={\textbf{Implementation of the Breadth-First Search (BFS) logic for the systematic downward construction of the solution landscape.}}, label=LandscapeRun.py(3), language=Python, showstringspaces=false]
    ... # continued from the preceding codes
	while tempList:
		tempsearchinitial = tempList.pop(0)
		for j in range(tempsearchinitial[2] - 1,
				max(0, tempsearchinitial[2] - instance.MaxIndexGap) - 1, -1):
			if j > 0:
				# Try all combinations of eigenvectors
				all_combinations = list(combinations(range(tempsearchinitial[2]), j))
				for r in range(len(all_combinations)):
					if instance.EigenCombination == "min":
						if r != 0:
							continue
					instance.calHiSD.SaddleIndex = j
					instance.calHiSD.InitialEigenVectors = tempsearchinitial[3][:, all_combinations[r]]
					# Generate perturbations
					PerturbationList = []
					for p in range(instance.PerturbationNumber):
    					pertemp = instance.PerMethod(instance.Dim, instance.PerturbationRadius)
    					pertemp_projection=tempsearchinitial[3]@(tempsearchinitial[3].T@pertemp)
    					pertemp_projection_normalized=instance.PerturbationRadius*pertemp_projection/np.clip(np.linalg.norm(pertemp_projection),1e-10,None)
    					PerturbationList.append(pertemp_projection_normalized)
    					PerturbationList.append(-pertemp_projection_normalized)
					for per in PerturbationList:
						instance.calHiSD.InitialPoint = tempsearchinitial[1] + per
						... # Formatted output
						instance.calHiSD.run()
						if not instance.calHiSD.flag:
							continue
						whetherinlist, smallestind = checkwhetherexist(instance, 
								instance.calHiSD.x, instance.calHiSD.finalindex)
						if instance.calHiSD.finalindex > j:
							print("Sorry! Because of the relaxed abort criteria, we find a saddle point with higher index.")
							continue
						if whetherinlist:  # If the point is already in the list
							xold = instance.SaddleList[smallestind][1]
							xnew = instance.calHiSD.x
							gold = instance.calHiSD.Grad(xold)
							gnew = instance.calHiSD.Grad(xnew)
							if np.linalg.norm(gold) > np.linalg.norm(gnew):
								instance.SaddleList[smallestind][1] = instance.calHiSD.x
							if (tempsearchinitial[0] not in instance.SaddleList[smallestind][4]):
								instance.SaddleList[smallestind][4].append(tempsearchinitial[0])
								if instance.SaveTrajectory:
									instance.DetailRecord.append([smallestind, tempsearchinitial[0],
											instance.calHiSD.x_record, instance.calHiSD.timestep_record])
								else:
									instance.DetailRecord.append([smallestind, tempsearchinitial[0]])
						else:  # Add the new saddle point
							instance.SaddleList.append([instance.saddleind, instance.calHiSD.x,
									instance.calHiSD.finalindex, instance.calHiSD.negvecs, [tempsearchinitial[0]]])
							if instance.SaveTrajectory:
								instance.DetailRecord.append([instance.saddleind, tempsearchinitial[0],
										instance.calHiSD.x_record, instance.calHiSD.timestep_record])
							else:
								instance.DetailRecord.append([instance.saddleind, tempsearchinitial[0]])
							tempList.append([instance.saddleind, instance.calHiSD.x,
									instance.calHiSD.finalindex, instance.calHiSD.negvecs])
							instance.saddleind = instance.saddleind + 1
			else:
				... # for index-0, no eigenvectors, perturbation only
\end{lstlisting}

\noindent\textbf{Perturbation}

As mentioned above, for downward searching, we invoke \graycode{Perturbation.py} to implement perturbation methods at saddle points (necessary for escaping current saddle points when searching for new ones from existing saddle points).
\begin{lstlisting}[caption={\textbf{Generation of normalized random perturbations to initiate escape from critical points along unstable directions.}}, label=Perturbation.py, language=Python]
def gaussianper(Dim, PerturbationRadius):
	"""
	Generate Gaussian-distributed random vector and normalize it.
	"""
	v = np.random.randn(Dim)
	v_normalized = PerturbationRadius * v / np.clip(np.linalg.norm(v), 1e-10, None)
	return v_normalized.reshape(-1, 1)


def uniformper(Dim, PerturbationRadius):
	"""
	Generate uniform-distributed random vector and normalize it.
	"""
	v = np.random.uniform(-1, 1, Dim)
	v_normalized = PerturbationRadius * v / np.clip(np.linalg.norm(v), 1e-10, None)
	return v_normalized.reshape(-1, 1)

... # Called in the core loop.
pertemp = instance.PerMethod(instance.Dim, instance.PerturbationRadius)
pertemp_projection = tempsearchinitial[3] @ (tempsearchinitial[3].T @ pertemp)
pertemp_projection_normalized = instance.PerturbationRadius * pertemp_projection / np.clip(np.linalg.norm(pertemp_projection), 1e-10, None)
PerturbationList.append(pertemp_projection_normalized)
PerturbationList.append(-pertemp_projection_normalized)
...
\end{lstlisting}

\noindent\textbf{Restart with Upward Search}

We can construct the complete solution landscape using upward search. There are two choices:
\begin{itemize}
  \item The restart function \graycode{RestartFromPoint} for re-solving with new initial points.
  \item The function \graycode{RestartFromSaddle}, which continues solving from an existing saddle point in the solution landscape (note that a small random perturbation is needed when restarting to escape from the saddle point).
\end{itemize}

\begin{lstlisting}[caption={\textbf{Methods for restarting the landscape construction from a specific point or an existing saddle.}}, label=Restart.py, language=Python]
class Landscape:
	... # other methods

	def RestartFromPoint(self, RestartPoint, MaxIndex):
		"""
		Restart the landscape process.
		"""
		RestartPoint = np.array(RestartPoint).reshape(-1, 1)
		self.calHiSD.InitialPoint = RestartPoint
		self.InitialPoint = RestartPoint
		historyMaxIndex = self.MaxIndex
		self.MaxIndex = MaxIndex
		self.BeginID = -1
		self.Continue = True
		LandscapeRun(self)
		self.Continue = False
		self.MaxIndex = historyMaxIndex

	def RestartFromSaddle(self, BeginID, Perturbation, MaxIndex):
		"""
		Run the landscape process.
		"""
		if BeginID >= len(self.SaddleList):
			raise ValueError("Invalid saddle ID")
		RestartPoint = self.SaddleList[BeginID][1] + np.array(Perturbation).reshape(
			-1, 1
		)
		self.calHiSD.InitialPoint = RestartPoint
		self.InitialPoint = RestartPoint
		historyMaxIndex = self.MaxIndex
		self.MaxIndex = MaxIndex
		self.BeginID = BeginID
		self.Continue = True
		LandscapeRun(self)
		self.Continue = False
		self.MaxIndex = historyMaxIndex
\end{lstlisting}

\subsection{Data Processing Module}
This module handles the saving and visualization of search results. As it primarily involves post-processing rather than numerical computation, we omit the code details and focus on the main workflow. Implementation details can be found in the GitHub repository.

\subsubsection{Search Trajectory Plotting Program \graycode{Trajectory.py}}
\label{Search Trajectory Plotting Program}
This program visualizes search trajectories with different approaches for varying data dimensions:
\begin{itemize}
  \item High-dimensional data: Use user-defined projection to 2D.
  \item 2D data: Generate grid points to compute energy values, plot contour lines (contour/contourf), mark saddle points with colors according to their indices, and draw search trajectories (detailed or straight-line connections).
  \item 1D data: Plot function curves, mark saddle point positions, and connect relevant saddle points with straight lines.
\end{itemize}

\subsubsection{Solution Landscape Plotting Program \graycode{Connection.py}}
\label{Solution Landscape Plotting Program}
This program generates topological connection diagrams of solution landscapes through these steps:
\begin{enumerate}
  \item Configure plot parameters.
  \item Construct graph structure.
  \item Calculate node coordinates (for plotting).
  \item Draw elements: Nodes from high to low Morse index (ensuring higher nodes appear above), dynamic node sizing, and edges.
  \item Add node labels and adjustments.
\end{enumerate}

\subsubsection{Data Saving Program \graycode{Save.py}}
\label{Data Saving Program}
This program saves the SaddleList data from instance objects in specified formats, including JSON (.json), MATLAB (.mat), and Python Pickle (.pkl) formats. It extracts each saddle point's ID, position, Morse index, and parent nodes. The serialization is dispatched according to the user-selected output format.

\section{Gallery}

In this section, we present four examples that demonstrate the capabilities of our implementation.

\subsection{Butterfly Function}
This example features a system whose energy heatmap resembles a butterfly. The energy function is given by
$$
\begin{aligned}
E(x,y)=x^{4}-1.5x^{2}y^{2}+y^{4}-2y^{3}+y^{2}+x^{2}y-2x^{2}.
\end{aligned}
$$
First, we initialize the package by importing the \texttt{Landscape} class.
\begin{lstlisting}[caption={\textbf{Importing necessary libraries and the main \texttt{Landscape} class.}}, label=Example-Butterfly.py(1), language=Python]
import sys
import os
sys.path.append(os.path.abspath(os.path.join(os.getcwd(), '..', 'saddlescape-1.0')))
from saddlescape import Landscape
import numpy as np
# import packages needed
\end{lstlisting}
We then define the symbolic energy function and initialize the parameters.
\begin{lstlisting}[caption={\textbf{Definition of the symbolic energy function and initialization of system parameters.}}, label=Example-Butterfly.py(2), language=Python]
# given energy function
energyfunction = 'x1**4 -1.5*x1**2*x2**2+ x2**4 - 2*x2**3 + x2**2 + x1**2*x2 - 2*x1**2' 
# parameter initialization
x0 = np.array([0.1, 0.1]) # initial point
dt = 1e-2 # time step
k = 2 # the maximum index of saddle point
acceme = 'none'
maxiter = 10000 # max iter
\end{lstlisting}
Next, we instantiate the solver and execute the \graycode{Run} method to construct the solution landscape.
\begin{lstlisting}[caption={\textbf{Instantiation of the solver and execution of the solution landscape construction.}}, label=Example-Butterfly.py(3), language=Python]
MyLandscape = Landscape(MaxIndex=k, AutoDiff=True, EnergyFunction=energyfunction, EigenMethod='euler', EigenMaxIter=1, InitialPoint=x0, TimeStep=dt, Acceleration=acceme, MaxIter=maxiter, EigenCombination='all', PerturbationNumber=1, PerturbationRadius=1e-2, Verbose=True, ReportInterval=100)
# Instantiation
MyLandscape.Run()
# Calculate
\end{lstlisting}
Partial output is shown below:
\begin{lstlisting}[caption={\textbf{Representative console output demonstrating initialization and iteration progress.}}, label=Output-Butterfly.py(1), language=Python]
  HiSD Solver Configuration:
  ------------------------------
  [HiSD] Current parameters (initialized):
  [Config Sync] `Dim` parameter auto-adjusted to 2 based on `InitialPoint` dimensionality.
  Parameter `NumericalGrad` not specified - using default value False.
  Using `EnergyFunction` instead of `Grad` - enabling auto-differentiation mode.
  Parameter `Momentum` not specified - using default value 0.0.
  Parameter `BBStep` not specified - using default value False.
  Parameter `DimerLength` not specified - using default value 1e-05.
  Parameter `Tolerance` not specified - using default value 1e-06.
  Parameter `NesterovChoice` not specified - using default value 1.
  Parameter `SearchArea` not specified - using default value 1000.0.
  
  ... # more parameter checks
  
  Gradient system detected. Activating HiSD algorithm.
  
  Landscape Configuration:
  ------------------------------
  [Landscape] Current parameters (initialized):
  Parameter `SameJudgementMethod` not specified - using default value <function LandscapeCheckParam.<locals>.<lambda> at 0x000002EBD1A2F240>.
  Parameter `PerturbationMethod` not specified - using default value uniform.
  Parameter `InitialEigenVectors` not specified - using default value None.
  Parameter `SaveTrajectory` not specified - using default value True.
  Parameter `MaxIndexGap` not specified - using default value 1.
  
  Start running:
  ------------------------------

  From initial point search index-2:
  ------------------------------
  
  Iteration: 100|| Norm of gradient: 0.176823
  Iteration: 200|| Norm of gradient: 0.084819
  ...
  Iteration: 1200|| Norm of gradient: 0.000004
  Iteration: 1300|| Norm of gradient: 0.000001
  Non-degenerate saddle point identified: Morse index =2 (number of negative eigenvalues).
  
  From saddle point (index-2, ID-0) search index-1:
  ------------------------------
  
  Iteration: 100|| Norm of gradient: 0.023111
  Iteration: 200|| Norm of gradient: 0.061258
  ...
  Iteration: 900|| Norm of gradient: 0.000024
  Iteration: 1000|| Norm of gradient: 0.000003
  Non-degenerate saddle point identified: Morse index =1 (number of negative eigenvalues).
  
  ... # more log messages during searches
\end{lstlisting}
After the searches, the data can be used for various kinds of post-processing. For instance, we draw the search trajectory.
\begin{lstlisting}[caption={\textbf{Command for visualizing the search trajectory.}}, label=Example-Butterfly.py(4), language=Python]
MyLandscape.DrawTrajectory(ContourGridNum=100, ContourGridOut=25, DetailedTraj=True)
# Draw the search path.
\end{lstlisting}
Also, drawing the solution landscape and saving the data are supported.
\begin{lstlisting}[caption={\textbf{Commands for visualizing the connectivity graph and saving results.}}, label=Example-Butterfly.py(5), language=Python]
MyLandscape.DrawConnection()
MyLandscape.Save('output/Ex_Butterfly')
# Save the data
\end{lstlisting}

\begin{figure}[H]
    \centering
    % 第一个子图 (a)
    \begin{minipage}[b]{0.45\linewidth}
        \centering
        % grid 参数用于调试位置，确认位置后删除 grid
        \begin{overpic}[width=\linewidth]{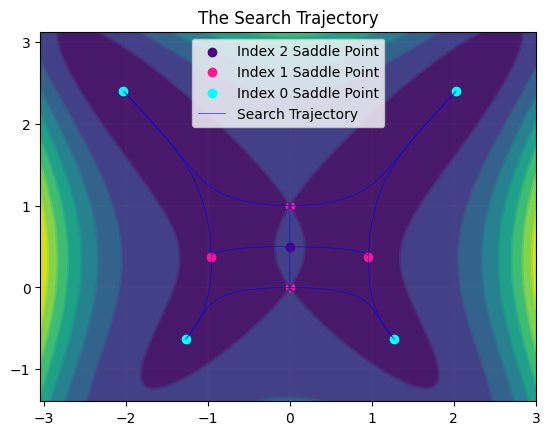}
            % put(x,y) 中的 x,y 是相对于图片宽高的百分比 (0-100)
            % \bfseries 加粗，\large 调整大小，根据需要调整
            \put(5,75){\bfseries\color{black}(a)} 
        \end{overpic}
        % 如果你不想用 \subfigure 的自动编号，可以手动写 caption，或者保留 \subfigure 仅作容器
        \vspace{0pt} % 图片和标题的间距

    \end{minipage}
    \hfill % 图片之间的间距
    % 第二个子图 (b)
    \begin{minipage}[b]{0.45\linewidth}
        \centering
        \begin{overpic}[width=\linewidth]{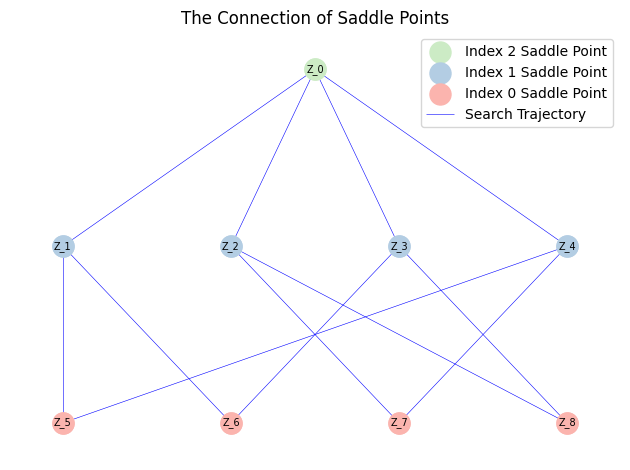}
            \put(5,75){\bfseries\color{black}(b)}
        \end{overpic}
        \vspace{0pt}

    \end{minipage}
    
    \caption{\textbf{Solution landscape reconstruction for the Butterfly potential energy surface.} \textbf{(a)} 2D projection of search trajectories initiated from $\boldsymbol{x}_0=(0.1, 0.1)$. \textbf{(b)} The solution landscape illustrating the connectivity between the index-2 (maximum), index-1 (transition states), and index-0 (minima) saddles.}
    \label{Figure-Butterfly(1)}
\end{figure}

Moreover, an animated image of the search can also be generated through specific post-processing. The details can be found at \\
\href{https://github.com/HiSDpackage/saddlescape/blob/main/_Examples/Ex_1_Butterfly.ipynb}{https://github.com/HiSDpackage/saddlescape/blob/main/\_Examples/Ex\_1\_Butterfly.ipynb}. 

The final output is provided as supplementary material and can be accessed at \\
\href{https://github.com/HiSDpackage/saddlescape/blob/main/_Examples/output/Ex_Butterfly.mp4}{https://github.com/HiSDpackage/saddlescape/blob/main/\_Examples/output/Ex\_Butterfly.mp4}.

\subsection{M\"uller-Brown Potential}

In this example, we demonstrate how to restart the algorithm to obtain a complete solution landscape via an upward search.

We define the M\"uller-Brown potential as a linear combination of Gaussian potentials:
$$
\begin{aligned}
E_{MB}(x,y)=\sum_{i=1}^{4}A_{i}\exp [a_{i}(x-\bar{x}_{i})^{2}+b_{i}(x-\bar{x}_{i})(y-\bar{y}_{i})+c_{i}(y-\bar{y}_{i})^{2}].
\end{aligned}
$$
The parameters are set as follows:
$$
\begin{aligned}
A &= [-200,-100,-170,15], \\
a &= [-1,-1,-6.5,0.7],\  b=[0,0,11,0.6],\  c=[-10,-10,-6.5,0.7], \\
\bar{x} &= [1,0,-0.5,-1],\  \bar{y}=[0,0.5,1.5,1].
\end{aligned}
$$

Similar to the first example, we initialize the solver (here starting from $(0.15, 0.25)$), search for saddle points, and then draw the search trajectory and the solution landscape.

\begin{figure}[H]
    \centering
    % 第一个子图 (a)
    \begin{minipage}[b]{0.45\linewidth}
        \centering
        % grid 参数用于调试位置，确认位置后删除 grid
        \begin{overpic}[width=\linewidth]{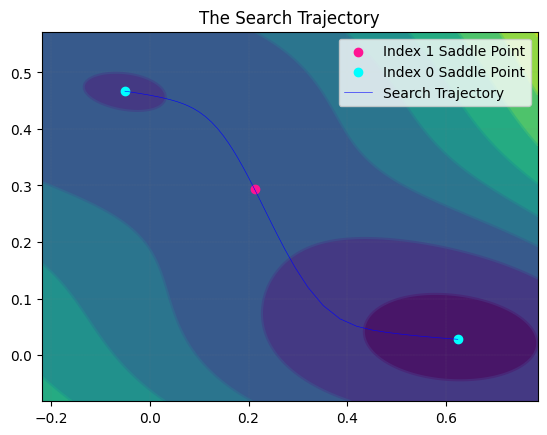}
            % put(x,y) 中的 x,y 是相对于图片宽高的百分比 (0-100)
            % \bfseries 加粗，\large 调整大小，根据需要调整
            \put(5,75){\bfseries\color{black}(a)} 
        \end{overpic}
        % 如果你不想用 \subfigure 的自动编号，可以手动写 caption，或者保留 \subfigure 仅作容器
        \vspace{0pt} % 图片和标题的间距

    \end{minipage}
    \hfill % 图片之间的间距
    % 第二个子图 (b)
    \begin{minipage}[b]{0.45\linewidth}
        \centering
        \begin{overpic}[width=\linewidth]{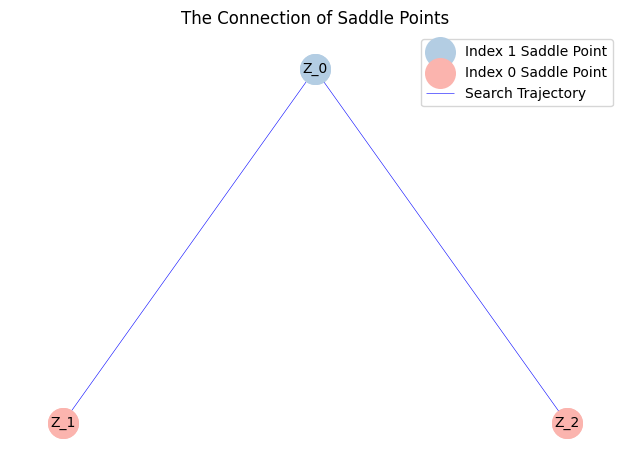}
            \put(5,75){\bfseries\color{black}(b)}
        \end{overpic}
        \vspace{0pt}

    \end{minipage}
    
    \caption{\textbf{Preliminary exploration of the M\"uller-Brown potential landscape.} The graphs \textbf{(a)} and \textbf{(b)} show an incomplete set of critical points and connections obtained from a single initial search pass, highlighting the need for iterative restart.}
    \label{Figure-Muller(1)}
\end{figure}

However, the M\"uller-Brown potential features a multimodal distribution, making the initial search incomplete. Therefore, we restart the search from a local minimum:
\begin{lstlisting}[caption={\textbf{Restarting the search from a specific saddle point to expand the landscape via upward search.}}, label=Example-Muller.py(1), language=Python]
MyLandscape.RestartFromSaddle(1,-np.array([[-0.01],[0]]),1)
# restart the search using an upward search
\end{lstlisting}
Also, we draw figures to analyze and save the data.
\begin{lstlisting}[caption={\textbf{Post-processing commands for visualization and data storage.}}, label=Example-Muller.py(2), language=Python]
MyLandscape.DrawTrajectory(ContourGridNum=100, ContourGridOut=25, DetailedTraj=True)
# Draw the search path. 
MyLandscape.DrawConnection()
MyLandscape.Save('output/Ex_MBP','pickle')
# Save the data
\end{lstlisting}
Now, we obtain a complete solution landscape.
\begin{figure}[H]
    \centering
    % 第一个子图 (a)
    \begin{minipage}[b]{0.45\linewidth}
        \centering
        % grid 参数用于调试位置，确认位置后删除 grid
        \begin{overpic}[width=\linewidth]{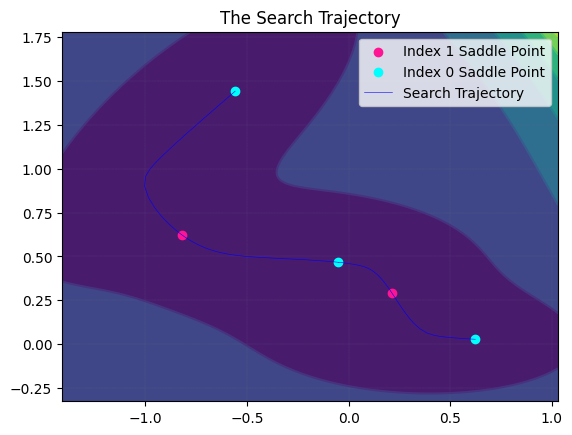}
            % put(x,y) 中的 x,y 是相对于图片宽高的百分比 (0-100)
            % \bfseries 加粗，\large 调整大小，根据需要调整
            \put(5,75){\bfseries\color{black}(a)} 
        \end{overpic}
        % 如果你不想用 \subfigure 的自动编号，可以手动写 caption，或者保留 \subfigure 仅作容器
        \vspace{0pt} % 图片和标题的间距

    \end{minipage}
    \hfill % 图片之间的间距
    % 第二个子图 (b)
    \begin{minipage}[b]{0.45\linewidth}
        \centering
        \begin{overpic}[width=\linewidth]{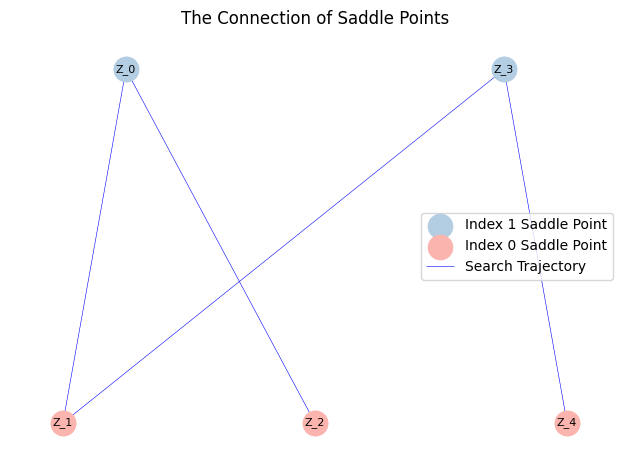}
            \put(5,75){\bfseries\color{black}(b)}
        \end{overpic}
        \vspace{0pt}

    \end{minipage}
    
    \caption{\textbf{The complete solution landscape of the M\"uller-Brown potential.} These comprehensive graphs \textbf{(a)} and \textbf{(b)} are obtained by restarting the search algorithm from the local minima identified in the preliminary phase (Figure \ref{Figure-Muller(1)}), successfully recovering the global topology.}
    \label{Figure-Muller(2)}
\end{figure}

\subsection{Cubic Function}

In this case, we demonstrate how to visualize search trajectories via projection, using a system with the energy function
$$
\begin{aligned}
E(x)=\sum_{j=1}^{n}j(x_{j}^{2}-1)^{2}.
\end{aligned}
$$
All saddle points are located at the vertices, face centers, and body center of an $n$-dimensional hypercube.

In particular, for $n=3$, there are 27 saddle points
$$
(a,b,c),\ \text{ where  }a,b,c \in \{-1,0,1\},
$$
which are the vertices, face centers, and body center of the cube.

We proceed as in the previous examples. Since the system is 3D, we define a projection function to visualize the trajectory.
\begin{lstlisting}[caption={\textbf{Definition of a projection function mapping 3D coordinates to 2D for visualization.}}, label=Example-Cubic.py(1), language=Python]
def proj_func(input):
  output = np.hstack((1.0 * input[:, [0]]+ 1.5 * input[:, [1]], 1.0 * input[:, [0]]+ 2.5 * input[:, [2]]))
  return output
\end{lstlisting}
Then we can draw the search trajectory via projection and the solution landscape.
\begin{lstlisting}[caption={\textbf{Commands to generate projected trajectory plots and the solution landscape for the cubic system.}}, label=Example-Cubic.py(2), language=Python]
MyLandscape.DrawTrajectory(ContourGridNum=100, ContourGridOut=25, Projection=proj_func)
# Draw the search path.
MyLandscape.DrawConnection()
MyLandscape.Save('output/Ex_Cubic','mat')
# Save the data
\end{lstlisting}

\begin{figure}[H]
    \centering
    % 第一个子图 (a)
    \begin{minipage}[b]{0.45\linewidth}
        \centering
        % grid 参数用于调试位置，确认位置后删除 grid
        \begin{overpic}[width=\linewidth]{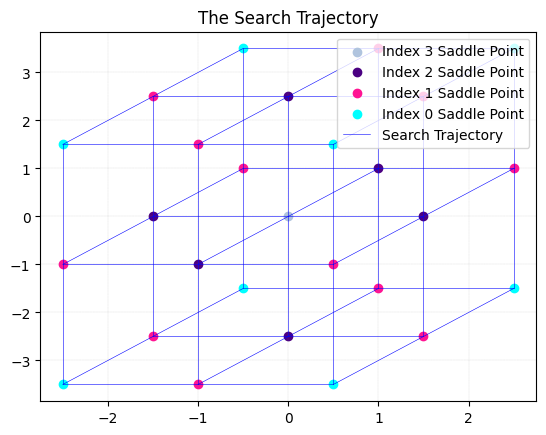}
            % put(x,y) 中的 x,y 是相对于图片宽高的百分比 (0-100)
            % \bfseries 加粗，\large 调整大小，根据需要调整
            \put(5,75){\bfseries\color{black}(a)} 
        \end{overpic}
        % 如果你不想用 \subfigure 的自动编号，可以手动写 caption，或者保留 \subfigure 仅作容器
        \vspace{0pt} % 图片和标题的间距

    \end{minipage}
    \hfill % 图片之间的间距
    % 第二个子图 (b)
    \begin{minipage}[b]{0.45\linewidth}
        \centering
        \begin{overpic}[width=\linewidth]{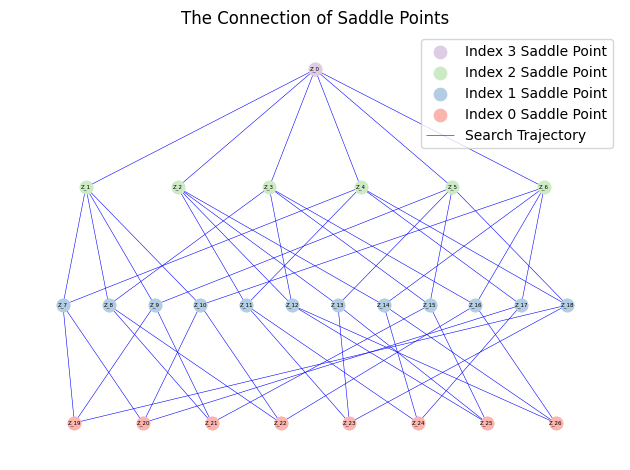}
            \put(5,75){\bfseries\color{black}(b)}
        \end{overpic}
        \vspace{0pt}

    \end{minipage}
    
    \caption{\textbf{Visualization of the 3D Cubic function landscape ($n=3$).} \textbf{(a)} Search trajectories projected onto a 2D plane using a linear combination of coordinates. \textbf{(b)} The fully connected solution landscape, exhibiting the symmetric structure characteristic of the hypercube vertices and face centers.}
    \label{Figure-Cubic(1)}
\end{figure}

\subsection{Phase Field Model}
In this example, we construct the solution landscape of a functional system. Specifically, we consider a 2D phase-field model:
$$
E(\phi) = \int_{\Omega} \left( \frac{\kappa}{2} |\nabla \phi|^2 + \frac{1}{4} (1 - \phi^2)^2 \right) dx.
$$
The Allen-Cahn equation is
$$
\dot{\phi} = \kappa \Delta \phi + \phi - \phi^3.
$$
We consider $\Omega = [0,1]^{2}$ with periodic boundary conditions and discretize the domain using a finite difference scheme on a $64 \times 64$ grid. The spatial derivatives are approximated using central finite differences implemented via convolution kernels.
\begin{lstlisting}[caption={\textbf{Implementation of the discrete gradient flow for the phase-field model using finite difference approximations.}}, label=Example-Phase.py(1), language=Python]
from scipy.ndimage import convolve

def gradient_PF(x, opt):
  kappa = opt['kappa']
  n2 = len(x)
  n = int(np.sqrt(n2))
  phi = x.reshape(n, n)
  h = opt['h']
  D2 = np.array([[0, 1, 0], [1, -4, 1], [0, 1, 0]]) / h**2
  conv_term = convolve(phi, D2, mode='wrap')
  F = -(kappa * conv_term + phi - phi**3)
  return F.reshape(n2, 1)
\end{lstlisting}
To account for periodic boundary conditions, we define an equivalence predicate with translational invariance.
\begin{lstlisting}[caption={\textbf{Saddle point equivalence check using a spatial translation approach for periodic boundary conditions.}}, label=Example-Phase.py(2), language=Python]
def SameSaddle(x, y):
	n= 64
	A = x.reshape(n,n)
	B = y.reshape(n,n)
	epsilon = 0.05
	row_indices = np.array([(np.arange(n) - dx) % n for dx in range(n)])
	col_indices = np.array([(np.arange(n) - dy) % n for dy in range(n)])
	row_shifted = np.zeros((n, n, n), dtype=B.dtype)
	for dx in range(n):
		row_shifted[dx] = B[row_indices[dx], :]
	for dx in range(n):
		for dy in range(n):
			shifted_B = row_shifted[dx][:, col_indices[dy]]
			max_diff = np.max(np.abs(A - shifted_B))
			if max_diff <= epsilon:
				return True
	return False
\end{lstlisting}
Alternatively, one can also use normalized cross-correlation computed via the Fast Fourier Transform (FFT) algorithm to define a translation-invariant equivalence predicate.
\begin{lstlisting}[caption={\textbf{Saddle point equivalence check using FFT algorithm to handle periodic shifts.}}, label=Example-Phase.py(3), language=Python]
from numpy.fft import fft2, ifft2

def SameSaddle(x, y):
    n = 64
    A = x.reshape(n, n)  
    B = y.reshape(n, n)  
    # 2D Fourier transform
    fft_A = fft2(A)
    fft_B = fft2(B)
    cross_corr_spectrum = np.conj(fft_A) * fft_B
    real_space_corr = np.real(ifft2(cross_corr_spectrum))
    max_corr = np.max(real_space_corr)
    norm_factor = np.sum(A**2)
    if norm_factor < 1e-7:
        return np.sum(B**2) < 1e-7
    normalized_max_corr = max_corr / norm_factor
    if normalized_max_corr > 0.99:
        return True
    return False
\end{lstlisting}
We focus on the system with $\kappa=0.02$.
\begin{lstlisting}[caption={\textbf{Configuration of phase-field system parameters and the gradient function.}}, label=Example-Phase.py(4), language=Python]
N = 64
opt = {'kappa': 0.02, 'N': N, 'h': 1/N,}
GradFunc = lambda x: gradient_PF(x, opt)
\end{lstlisting}
We initialize the solver and solve the system.
\begin{lstlisting}[caption={\textbf{Initialization and execution of the HiSD solver for the phase-field model.}}, label=Example-Phase.py(5), language=Python]
# parameter initialization
x0 = np.array([0 for i in range(N**2)]) # initial point
dt = 1e-3 # time step
k = 5
acceme = 'nesterov'
neschoice = 1
nesres = 200
mom = 0.8
maxiter = 2000 # max iter

MyLandscape = Landscape(MaxIndex=k, AutoDiff=False, Grad=GradFunc, DimerLength=1e-3, HessianDimerLength=1e-3, EigenStepSize=1e-7, InitialPoint=x0, TimeStep=dt, Acceleration=acceme, SearchArea=1e4, SymmetryCheck=False,Tolerance=1e-4, MaxIndexGap=3, EigenCombination='min',BBStep=True, NesterovChoice=neschoice, NesterovRestart=nesres, Momentum=mom, MaxIter=maxiter, Verbose=True, ReportInterval=10, EigenMaxIter=2, PerturbationNumber=2,EigvecUnified=True,SaveTrajectory=False,SameJudgementMethod=SameSaddle, PerturbationRadius=5.0)
# Instantiation
MyLandscape.Run()
# Calculate
\end{lstlisting}
We can draw the solution landscape and save the data.
\begin{lstlisting}[caption={\textbf{Commands for saving and plotting the phase-field landscape results.}}, label=Example-Phase.py(6), language=Python]
MyLandscape.DrawConnection()
MyLandscape.Save('output/Ex_PhaseField','mat')
# Save the data
\end{lstlisting}
 The saddle points, representing functions, can be visualized via post-processing:
\begin{lstlisting}[caption={\textbf{Helper function to visualize the spatial phase field configuration $\phi(\boldsymbol{x})$ as a heatmap.}}, label=Example-Phase.py(7), language=Python]
import matplotlib.pyplot as plt
from scipy.interpolate import griddata

def plot_phi_heatmap(phi_vector, N):
	if len(phi_vector) != N * N:
		raise ValueError(f"Input shape must be {N * N}, but got {len(phi_vector)}")
	h=1/N
	x = np.linspace(h/2, 1-h/2, N)
	y = np.linspace(h/2, 1-h/2, N)
	X, Y = np.meshgrid(x, y)
	phi = phi_vector.reshape((N, N))
	grid_x, grid_y = np.mgrid[0:1:N*10j, 0:1:N*10j] 
	grid_phi = griddata((X.flatten(), Y.flatten()), phi.flatten(), (grid_x, grid_y), method='cubic')
	# Draw heatmap
	plt.figure(figsize=(8, 6))
	plt.imshow(grid_phi, extent=(0, 1, 0, 1), origin='lower', cmap='viridis',vmin=-1, vmax=1)
	plt.colorbar(label='$\phi$ value')
	plt.title('phase field $\phi$ heatmap')
	plt.xlabel('x')
	plt.ylabel('y')
	plt.show()

for i in range(len(MyLandscape.SaddleList)):
	plot_phi_heatmap(MyLandscape.SaddleList[i][1],N)
\end{lstlisting}
By visualizing saddle‑point phase fields alongside the solution landscape, we illustrate their correspondence.

\begin{figure}[H]
    \centering
    % 第一个子图 (a)
    \begin{minipage}[b]{0.45\linewidth}
        \centering
        % grid 参数用于调试位置，确认位置后删除 grid
        \begin{overpic}[width=\linewidth]{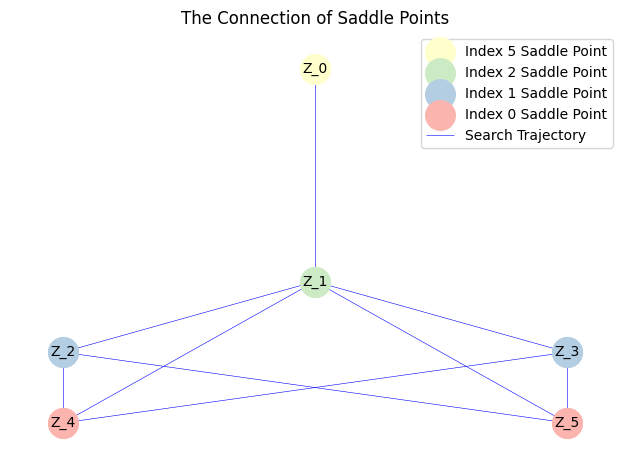}
            % put(x,y) 中的 x,y 是相对于图片宽高的百分比 (0-100)
            % \bfseries 加粗，\large 调整大小，根据需要调整
            \put(5,75){\bfseries\color{black}(a)} 
        \end{overpic}
        % 如果你不想用 \subfigure 的自动编号，可以手动写 caption，或者保留 \subfigure 仅作容器
        \vspace{0pt} % 图片和标题的间距

    \end{minipage}
    \hfill % 图片之间的间距
    % 第二个子图 (b)
    \begin{minipage}[b]{0.45\linewidth}
        \centering
        \begin{overpic}[width=\linewidth]{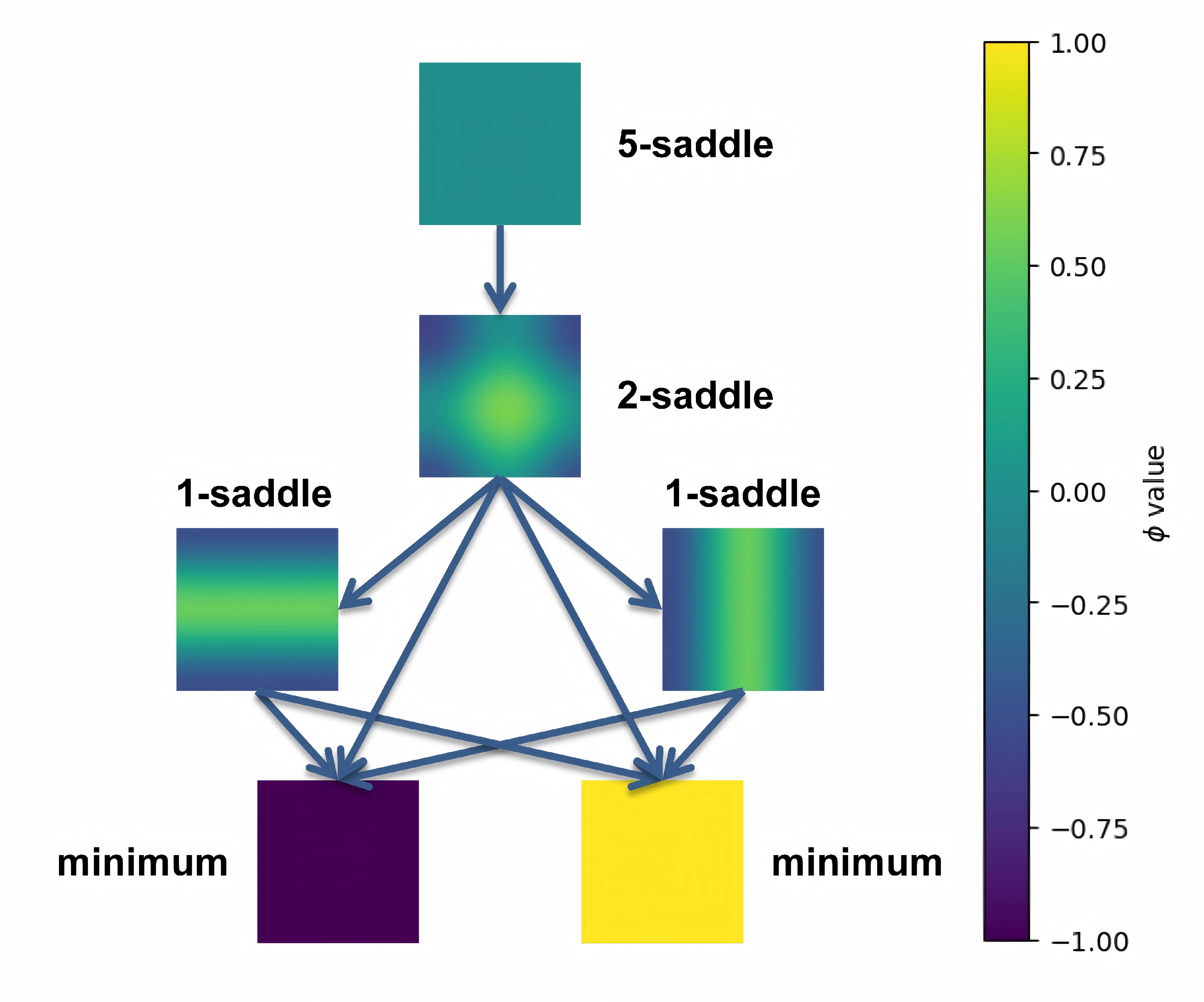}
            \put(5,75){\bfseries\color{black}(b)}
        \end{overpic}
        \vspace{0pt}

    \end{minipage}
    
    \caption{\textbf{Computed solution landscape for the phase-field model with interaction parameter $\kappa=0.02$.} \textbf{(a)} The solution landscape. \textbf{(b)} Inset images visualize the spatial phase configurations $\phi(\boldsymbol{x})$ corresponding to selected saddle points.}
    \label{Figure-Phase(1)}
\end{figure}

Moreover, for a system with $\kappa=0.01$ that possesses a more complex solution landscape, we can initially employ the fast BB step size combined with Nesterov acceleration to obtain an incomplete solution landscape with high efficiency. Subsequently, we restart the search using a fixed small step size to search locally and construct the complete solution landscape. This two-stage strategy, i.e., accelerated search followed by refined exploration, effectively balances computational efficiency and solution accuracy. The final solution landscape obtained through this approach is presented as follows:

\begin{figure}[H]
    \centering
    % 第一个子图 (a)
    \begin{minipage}[b]{0.45\linewidth}
        \centering
        % grid 参数用于调试位置，确认位置后删除 grid
        \begin{overpic}[width=\linewidth]{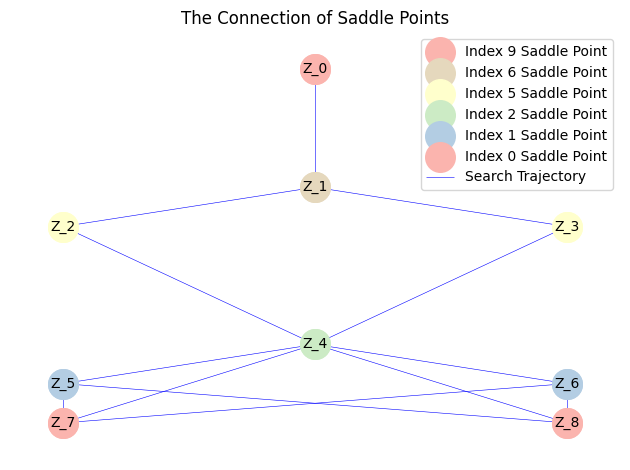}
            % put(x,y) 中的 x,y 是相对于图片宽高的百分比 (0-100)
            % \bfseries 加粗，\large 调整大小，根据需要调整
            \put(5,75){\bfseries\color{black}(a)} 
        \end{overpic}
        % 如果你不想用 \subfigure 的自动编号，可以手动写 caption，或者保留 \subfigure 仅作容器
        \vspace{0pt} % 图片和标题的间距

    \end{minipage}
    \hfill % 图片之间的间距
    % 第二个子图 (b)
    \begin{minipage}[b]{0.45\linewidth}
        \centering
        \begin{overpic}[width=\linewidth]{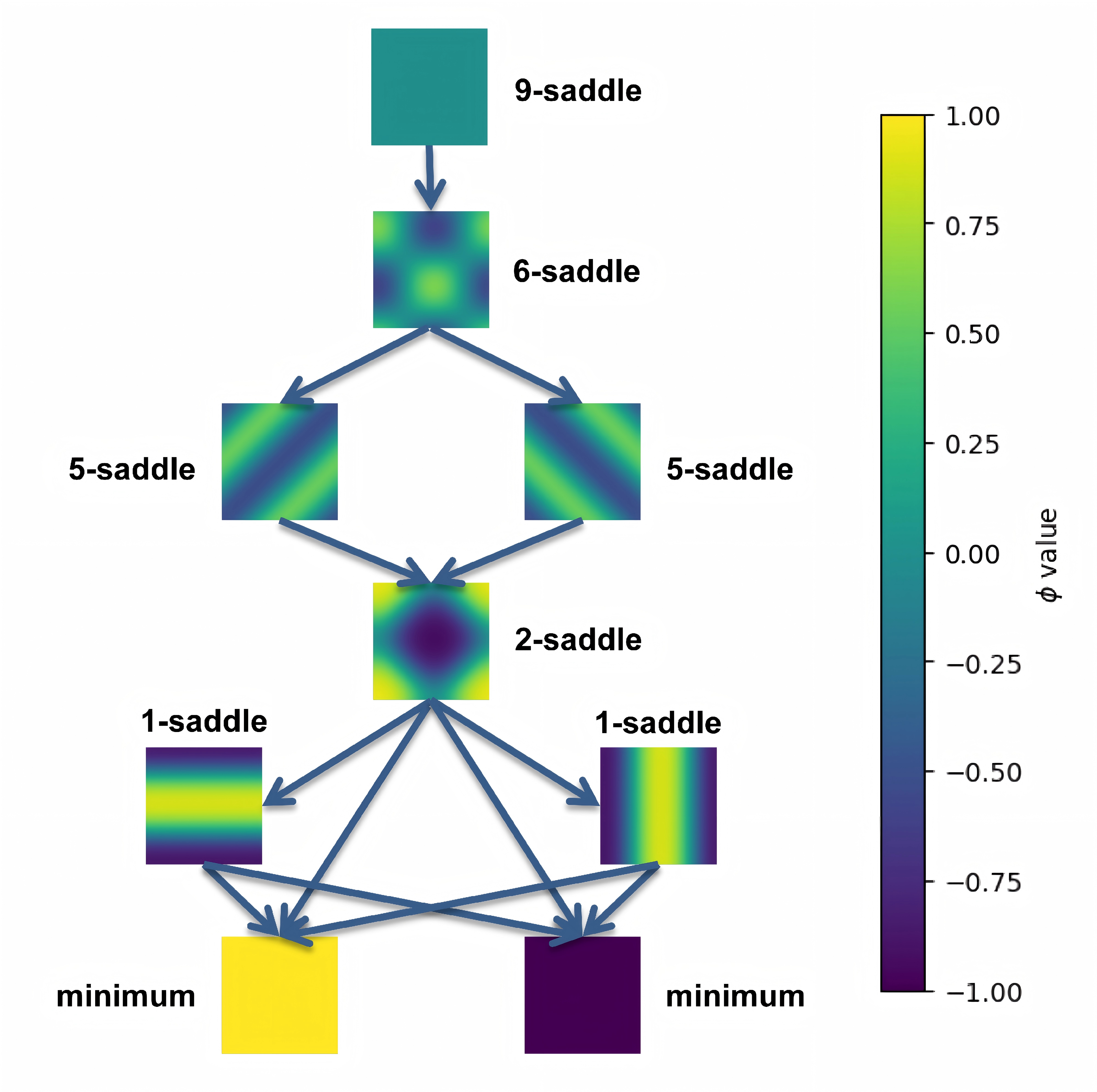}
            \put(5,75){\bfseries\color{black}(b)}
        \end{overpic}
        \vspace{0pt}

    \end{minipage}
    
    \caption{\textbf{Complexity increase in the solution landscape for the phase-field model at $\kappa=0.01$.} The reduced parameter value leads to a proliferation of metastable states and high-index saddles compared to the $\kappa=0.02$ case.}
    \label{Figure-Phase(2)}
\end{figure}

\section{Conclusion and Future Work}

\subsection{Conclusion}

The software package \texttt{SaddleScape V1.0} introduces a unified computational framework for exploring solution landscapes through advanced saddle dynamics algorithms, namely HiSD and its variants. Designed to balance accessibility and flexibility, the package enables both novice users and domain experts to analyze complex systems and construct solution landscapes with minimal effort. Key achievements of the package include:

\begin{itemize}
    \item \textbf{Automated Solution Landscape Construction:} A systematic, hierarchical, and largely automated procedure for constructing solution landscapes, based on downward and upward searches, randomized perturbations, and a BFS strategy with FIFO queue scheduling. This replaces the traditional practice of manually designing separate searches for each connection between saddle points.

    \item \textbf{Unified Framework and Generalized Input Support:} A unified implementation of HiSD and its variants (such as GHiSD for non-gradient systems and AHiSD with momentum acceleration) within a single interface. The package is compatible with both gradient and non-gradient systems via analytical expressions or Python callables, including functional systems after spatial discretization, with default parameters provided to reduce user intervention.

    \item \textbf{Flexible and Efficient Numerical Implementations:} Decoupled numerical modules (e.g., \graycode{AnaCal}, \graycode{HessianMatrix}, \graycode{EigMethod}) allow users to customize gradient evaluation, Hessian/Jacobian operators, eigenvalue solvers, and step-size or acceleration strategies. Symbolic computation, automatic differentiation, and matrix-free Hessian/Jacobian-vector products based on dimer-type finite differences are integrated with efficient vectorized linear-algebra routines to enhance robustness and computational efficiency, particularly for high-dimensional problems.

    \item \textbf{Landscape-level Control, Modularity, and Usability:} Users can configure system definitions, solver parameters, maximum index, perturbation strategies, and equivalence criteria for identifying distinct saddles, as well as restart procedures from existing saddles or new initial points, thereby refining or extending a partially constructed landscape. The package is organized into modular components for parameter checking, derivative handling, eigenvalue computation, iterative saddle dynamics solvers, and solution landscape construction. Default values are provided for many parameters, lowering the barrier for new users while still allowing advanced users to fine-tune numerical details and experiment with more advanced algorithms by replacing the corresponding modules.

    \item \textbf{Integrated Postprocessing, Visualization, and Data Export:} Automated workflows for search trajectory visualization, solution landscape construction, and data export in multiple formats (JSON, MATLAB, and Python Pickle) provide an intuitive understanding of both the search process and the resulting solution landscape, and facilitate downstream analysis and integration with other computational tools.
\end{itemize}

By abstracting numerical and algorithmic complexities into intuitive interfaces and automated workflows, \texttt{SaddleScape V1.0} bridges theoretical saddle dynamics with practical applications in areas such as condensed matter physics, materials science, biophysical modeling, and nonlinear dynamics, enabling researchers to focus on domain-specific insights rather than implementation details.

\subsection{Future Work}

To advance the \texttt{SaddleScape} package, several directions merit further exploration. On the algorithmic side, a natural extension is to integrate constrained variants such as CHiSD \cite{yin2022constrained}, data-driven formulations based on surrogate models \cite{JJIAM2023,Liu2024Neural}, and infinite-dimensional extensions such as Functional HiSD \cite{Zhang2025Numerical}. Another important direction is developing parallel HiSD algorithms by leveraging distributed computing frameworks to construct solution landscapes of high-dimensional complex systems more efficiently.

From a technical perspective, GPU-accelerated linear algebra and sparse matrix operations can be used for large-scale problems, and it is desirable to extend automatic differentiation support to stochastic and implicit dynamical systems. In terms of usability, future work includes the development of a graphical user interface (GUI) and tighter integration with Jupyter Notebooks, together with collaborative platforms for sharing pre-configured workflows and user-contributed modules.

Finally, there is substantial potential for interdisciplinary applications. Promising targets include materials science (e.g., phase transitions), biochemistry (e.g., protein folding), and machine learning (e.g., loss landscape analysis). In this context, inverse-design frameworks provide an important route for experimental validation, by first specifying desired macroscopic properties or target configurations and then using optimization or learning-based algorithms on the solution landscape to identify candidate system parameters that can subsequently be realized and tested in the laboratory.

The \texttt{SaddleScape} package will continue to evolve as an open-source platform, fostering collaboration between computational mathematics and applied sciences. Its modular design and community-driven development aim to establish it as a cornerstone toolkit for solution landscape exploration across disciplines.

\section{Acknowledgments}

We are grateful to Prof. Pingwen Zhang and Prof. Jianyuan Yin for their contributions to the HiSD algorithm. We would also like to acknowledge Dr. Bing Yu, Dr. Yue Luo, Prof. Xiangcheng Zheng, Dr. Baoming Shi, Yuankai Liu, and Haoran Wang for their significant contributions to the promotion and practical implementation of the algorithm.

\section{Funding}

This work was supported by the National Natural Science Foundation of China (No. 12225102, T2321001, and 12288101), the National Key Research and Development Program of China 2024YFA0919500, JR25003, and the Beijing Outstanding Young Scientist Program (No.JWZQ20240101027).

\section{Code Availability }

The source code can be found at the repository \href{https://github.com/HiSDpackage/saddlescape}{https://github.com/HiSDpackage/saddlescape}. The introductory website of the package is available at \href{https://hisdpackage.github.io/saddlescape}{https://hisdpackage.github.io/saddlescape}.

\bibliographystyle{siam}
\bibliography{main}

\begin{thebibliography}{10}

\bibitem{transition-rates1}
{\sc J.~Baker}, {\em An algorithm for the location of transition states},
  Journal of Computational Chemistry, 7 (1986), pp.~385--395.

\bibitem{Barzilai1988}
{\sc J.~Barzilai and J.~M. Borwein}, {\em Two-point step size gradient
  methods}, IMA Journal of Numerical Analysis, 8 (1988), pp.~141--148.

\bibitem{homotopy-methods2}
{\sc C.~Chen and Z.~Xie}, {\em Search extension method for multiple solutions
  of a nonlinear problem}, Computers \& Mathematics with Applications, 47
  (2004), pp.~327--343.

\bibitem{KineticEq}
{\sc H.~Chen, S.~Kandasamy, S.~Orszag, R.~Shock, S.~Succi, and V.~Yakhot}, {\em
  Extended boltzmann kinetic equation for turbulent flows}, Science, 301
  (2003), pp.~633--636.

\bibitem{BioChem1}
{\sc W.~W. Chen, M.~Niepel, and P.~K. Sorger}, {\em Classic and contemporary
  approaches to modeling biochemical reactions}, Genes \& Development, 24
  (2010), pp.~1861--1875.

\bibitem{PhysRevLett.104.148301}
{\sc X.~Cheng, L.~Lin, W.~E, P.~Zhang, and A.-C. Shi}, {\em Nucleation of
  ordered phases in block copolymers}, Physical Review Letters, 104 (2010),
  p.~148301.

\bibitem{neural-networks1}
{\sc R.~Das and D.~J. Wales}, {\em Energy landscapes for a machine-learning
  prediction of patient discharge}, Physical Review E, 93 (2016), p.~063310.

\bibitem{neural-networks2}
{\sc F.~Draxler, K.~Veschgini, M.~Salmhofer, and F.~Hamprecht}, {\em
  Essentially no barriers in neural network energy landscape}, in Proceedings
  of the 35th International Conference on Machine Learning, vol.~80 of PMLR,
  Jul 2018, pp.~1309--1318.

\bibitem{neural-networks3}
{\sc W.~E, C.~Ma, and L.~Wu}, {\em A comparative analysis of optimization and
  generalization properties of two-layer neural network and random feature
  models under gradient descent dynamics}, Science China Mathematics, 63
  (2020), pp.~1235--1258.

\bibitem{deflation-techniques}
{\sc P.~E. Farrell, A.~Birkisson, and S.~W. Funke}, {\em Deflation techniques
  for finding distinct solutions of nonlinear partial differential equations},
  SIAM Journal on Scientific Computing, 37 (2015), pp.~A2026--A2045.

\bibitem{han2022elastic}
{\sc Y.~Han, J.~Harris, A.~Majumdar, and L.~Zhang}, {\em Elastic anisotropy in
  the reduced landau-de gennes model}, Proceedings of the Royal Society A:
  Mathematical, Physical and Engineering Sciences, 478 (2022), p.~20210966.

\bibitem{han2021a}
{\sc Y.~Han, J.~Yin, Y.~Hu, A.~Majumdar, and L.~Zhang}, {\em Solution
  landscapes of the simplified ericksen-leslie model and its comparison with
  the reduced landau-de gennes model}, Proceedings of the Royal Society
  A-mathematical Physical and Engineering Science, 477 (2021), p.~20210458.

\bibitem{han2021solution}
{\sc Y.~Han, J.~Yin, P.~Zhang, A.~Majumdar, and L.~Zhang}, {\em Solution
  landscape of a reduced {L}andau--de {G}ennes model on a hexagon},
  Nonlinearity, 34 (2021), pp.~2048--2069.

\bibitem{transition-rates2}
{\sc A.~Heyden, A.~T. Bell, and F.~J. Keil}, {\em Efficient methods for finding
  transition states in chemical reactions: comparison of improved dimer method
  and partitioned rational function optimization method}, The Journal of
  Chemical Physics, 123 (2005), p.~224101.

\bibitem{protein-folding1}
{\sc D.~T. Leeson, F.~Gai, H.~M. Rodriguez, L.~M. Gregoret, and R.~B. Dyer},
  {\em Protein folding and unfolding on a complex energy landscape},
  Proceedings of the National Academy of Sciences, 97 (2000), pp.~2527--2532.

\bibitem{li2024acta}
{\sc L.~Li, B.~Yu, P.~Gao, J.~Lv, L.~Zhang, Y.~Wang, and Y.~Ma}, {\em
  Representing crystal potential energy surfaces via a stationary-point
  network}, Acta Materialia, 281 (2024), p.~120403.

\bibitem{Liu2024Neural}
{\sc Y.~Liu, L.~Zhang, and J.~Zhao}, {\em Neural network-based high-index
  saddle dynamics method for searching saddle points and solution landscape}.
\newblock arXiv:2411.16200, 2024.

\bibitem{luo2024NMPDE}
{\sc Y.~Luo, L.~Zhang, P.~Zhang, Z.~Zhang, and X.~Zheng}, {\em Semi-implicit
  method of high-index saddle dynamics and application to construct solution
  landscape}, Numerical Methods for Partial Differential Equations, 40 (2024),
  p.~e23123.

\bibitem{luo2025}
{\sc Y.~Luo, L.~Zhang, and X.~Zheng}, {\em Accelerated high-index saddle
  dynamics method for searching high-index saddle points}, Journal of
  Scientific Computing, 102 (2025), p.~31.

\bibitem{luo2022sinum}
{\sc Y.~Luo, X.~Zheng, X.~Cheng, and L.~Zhang}, {\em Convergence analysis of
  discrete high-index saddle dynamics}, SIAM Journal on Numerical Analysis, 60
  (2022), pp.~2731--2750.

\bibitem{protein-folding2}
{\sc F.~Mallamace, C.~Corsaro, D.~Mallamace, S.~Vasi, C.~Vasi, P.~Baglioni,
  S.~V. Buldyrev, S.-H. Chen, and H.~E. Stanley}, {\em Energy landscape in
  protein folding and unfolding}, Proceedings of the National Academy of
  Sciences, 113 (2016), pp.~3159--3163.

\bibitem{homotopy-methods1}
{\sc D.~Mehta}, {\em Finding all the stationary points of a potential-energy
  landscape via numerical polynomial-homotopy-continuation method}, Physical
  Review E, 84 (2011), p.~025702.

\bibitem{miao2025}
{\sc S.~Miao, L.~Zhang, P.~Zhang, and X.~Zheng}, {\em Construction and analysis
  for orthonormalized runge–kutta schemes of high-index saddle dynamics},
  Communications in Nonlinear Science and Numerical Simulation, 145 (2025),
  p.~108731.

\bibitem{Milnor+1963}
{\sc J.~Milnor}, {\em Morse Theory. (AM-51), Volume 51}, Princeton University
  Press, Princeton, 1963.

\bibitem{BioChem3}
{\sc Q.~Nie, L.~Qiao, Y.~Qiu, L.~Zhang, and W.~Zhao}, {\em Noise control and
  utility: From regulatory network to spatial patterning}, Science China
  Mathematics, 63 (2020), pp.~425--440.

\bibitem{protein-folding3}
{\sc J.~N. Onuchic, Z.~Luthey-Schulten, and P.~G. Wolynes}, {\em Theory of
  protein folding: The energy landscape perspective}, Annual Review of Physical
  Chemistry, 48 (1997), pp.~545--600.

\bibitem{BioChem2}
{\sc L.~Qiao, W.~Zhao, C.~Tang, Q.~Nie, and L.~Zhang}, {\em Network topologies
  that can achieve dual function of adaptation and noise attenuation}, Cell
  Systems, 9 (2019), pp.~271--285.e7.

\bibitem{shi2025siap}
{\sc B.~Shi, Y.~Han, C.~Ma, A.~Majumdar, and L.~Zhang}, {\em A modified
  landau–de gennes theory for smectic liquid crystals: Phase transitions and
  structural transitions}, SIAM Journal on Applied Mathematics, 85 (2025),
  pp.~821--847.

\bibitem{shi2024}
{\sc B.~Shi, Y.~Han, A.~Majumdar, and L.~Zhang}, {\em Multistability for
  nematic liquid crystals in cuboids with degenerate planar boundary
  conditions}, SIAM Journal on Applied Mathematics, 84 (2024), pp.~756--781.

\bibitem{shi2023nonlinearity}
{\sc B.~Shi, Y.~Han, J.~Yin, A.~Majumdar, and L.~Zhang}, {\em Hierarchies of
  critical points of a landau-de gennes free energy on three-dimensional
  cuboids}, Nonlinearity, 36 (2023), p.~2631.

\bibitem{shi2022siap}
{\sc B.~Shi, Y.~Han, and L.~Zhang}, {\em Nematic liquid crystals in a
  rectangular confinement: Solution landscape, and bifurcation}, SIAM Journal
  on Applied Mathematics, 82 (2022), pp.~1808--1828.

\bibitem{miao2025csiam}
{\sc M.~Shuai, L.~Ziqi, L.~Zhang, P.~Zhang, and X.~Zheng}, {\em Construction
  and analysis for adams explicit discretization of high-index saddle
  dynamics}, CSIAM Transactions on Applied Mathematics, 6 (2025), pp.~625--650.

\bibitem{su2025sinum}
{\sc H.~Su, H.~Wang, L.~Zhang, J.~Zhao, and X.~Zheng}, {\em Improved high-index
  saddle dynamics for finding saddle points and solution landscape}, SIAM
  Journal on Numerical Analysis, 63 (2025), pp.~1757--1775.

\bibitem{NavierStokes}
{\sc R.~Temam}, {\em Navier-Stokes Equations: Theory and Numerical Analysis},
  vol.~2 of Studies in Mathematics and its Applications, Elsevier Science
  Publishers B.V., Amsterdam, 3rd revised~ed., 1984.

\bibitem{biological-applications1}
{\sc D.~Wales}, {\em Energy Landscapes: Applications to Clusters, Biomolecules
  and Glasses}, Cambridge Molecular Science, Cambridge University Press,
  Cambridge, 2004.

\bibitem{atomic}
{\sc D.~J. Wales and J.~P.~K. Doye}, {\em Global optimization by basin-hopping
  and the lowest energy structures of lennard-jones clusters containing up to
  110 atoms}, The Journal of Physical Chemistry A, 101 (1997), pp.~5111--5116.

\bibitem{wang2025siap}
{\sc S.~Wang, T.~Wang, S.~Wu, L.~Zhang, and X.~Zou}, {\em Mathematical modeling
  and solution landscape reveal cancer progression dynamics in tumor ecological
  microenvironment}, SIAM Journal on Applied Mathematics, 85 (2025),
  pp.~50--77.

\bibitem{wang2021modeling}
{\sc W.~Wang, L.~Zhang, and P.~Zhang}, {\em Modelling and computation of liquid
  crystals}, Acta Numerica, 30 (2021), pp.~765--851.

\bibitem{xu2021solution}
{\sc Z.~Xu, Y.~Han, J.~Yin, B.~Yu, Y.~Nishiura, and L.~Zhang}, {\em Solution
  landscapes of the diblock copolymer-homopolymer model under two-dimensional
  confinement}, Physical Review E, 104 (2021), p.~014505.

\bibitem{xue2024elife}
{\sc G.~Xue, X.~Zhang, W.~Li, L.~Zhang, Z.~Zhang, X.~Zhou, D.~Zhang, L.~Zhang,
  and Z.~Li}, {\em A logic-incorporated gene regulatory network deciphers
  principles in cell fate decisions}, eLife, 12 (2024), p.~RP88742.

\bibitem{yao2022cicp}
{\sc X.~Yao, J.~Xu, and L.~Zhang}, {\em Transition pathways in cylinder-gyroid
  interface}, Communications in Computational Physics, 32 (2022), pp.~810--828.

\bibitem{yin2024revealing}
{\sc J.~Yin, Z.~Huang, Y.~Cai, Q.~Du, and L.~Zhang}, {\em Revealing excited
  states of rotational bose-einstein condensates}, The Innovation, 5 (2024).

\bibitem{yin2022constrained}
{\sc J.~Yin, Z.~Huang, and L.~Zhang}, {\em Constrained high-index saddle
  dynamics for the solution landscape with equality constraints}, Journal of
  Scientific Computing, 91 (2022), p.~62.

\bibitem{Yin2021transition}
{\sc J.~Yin, K.~Jiang, A.-C. Shi, P.~Zhang, and L.~Zhang}, {\em Transition
  pathways connecting crystals and quasicrystals}, Proceedings of the National
  Academy of the Sciences of the United States of America, 118 (2021),
  p.~e2106230118.

\bibitem{yin2020construction}
{\sc J.~Yin, Y.~Wang, J.~Z.~Y. Chen, P.~Zhang, and L.~Zhang}, {\em Construction
  of a pathway map on a complicated energy landscape}, Physical Review Letters,
  124 (2020), p.~090601.

\bibitem{yin2021searching}
{\sc J.~Yin, B.~Yu, and L.~Zhang}, {\em Searching the solution landscape by
  generalized high-index saddle dynamics}, Science China Mathematics, 64
  (2021), pp.~1801--1816.

\bibitem{yin2019high}
{\sc J.~Yin, L.~Zhang, and P.~Zhang}, {\em High-index optimization-based
  dhrinking dimer method for finding high-index saddle points}, SIAM Journal on
  Scientific Computing, 41 (2019), pp.~A3576--A3595.

\bibitem{yin2022solution}
\leavevmode\vrule height 2pt depth -1.6pt width 23pt, {\em Solution landscape
  of the onsager model identifies non-axisymmetric critical points}, Physica D,
  430 (2022), p.~133081.

\bibitem{yu2022jcp}
{\sc B.~Yu, X.~Zheng, P.~Zhang, and L.~Zhang}, {\em Computing solution
  landscape of nonlinear space-fractional problems via fast approximation
  algorithm}, Journal of Computational Physics, 468 (2022), p.~111513.

\bibitem{biological-applications2}
{\sc P.~Yu, Q.~Nie, C.~Tang, and L.~Zhang}, {\em Nanog induced intermediate
  state in regulating stem cell differentiation and reprogramming}, BMC Systems
  Biology, 12 (2018), p.~22.

\bibitem{jcm2023}
{\sc L.~Zhang}, {\em Construction of solution landscapes for complex systems},
  Mathematica Numerica Sinica, 45 (2023), pp.~267--283.

\bibitem{PhysRevLett.98.265703}
{\sc L.~Zhang, L.-Q. Chen, and Q.~Du}, {\em Morphology of critical nuclei in
  solid-state phase transformations}, Physical Review Letters, 98 (2007),
  p.~265703.

\bibitem{biological-applications3}
{\sc L.~Zhang, K.~Radtke, L.~Zheng, A.~Q. Cai, T.~F. Schilling, and Q.~Nie},
  {\em Noise drives sharpening of gene expression boundaries in the zebrafish
  hindbrain}, Molecular Systems Biology, 8 (2012), p.~613.

\bibitem{zhang2022sinum}
{\sc L.~Zhang, P.~Zhang, and X.~Zheng}, {\em Error estimates for euler
  discretization of high-index saddle dynamics}, SIAM Journal on Numerical
  Analysis, 60 (2022), pp.~2925--2944.

\bibitem{scm2023}
\leavevmode\vrule height 2pt depth -1.6pt width 23pt, {\em Discretization and
  index-robust error analysis for constrained high-index saddle dynamics on the
  high-dimensional sphere}, Science China Mathematics, 66 (2023),
  pp.~2347--2360.

\bibitem{cam2023}
\leavevmode\vrule height 2pt depth -1.6pt width 23pt, {\em Error estimate for
  semi-implicit method of sphere-constrained high-index saddle dynamics},
  Chinese Annals of Mathematics, Series B, 44 (2023), pp.~765--780.

\bibitem{csiam2023}
\leavevmode\vrule height 2pt depth -1.6pt width 23pt, {\em Mathematical and
  numerical analysis to shrinking-dimer saddle dynamics with local lipschitz
  conditions}, CSIAM Transactions on Applied Mathematics, 4 (2023),
  pp.~157--176.

\bibitem{JJIAM2023}
\leavevmode\vrule height 2pt depth -1.6pt width 23pt, {\em A model-free
  shrinking-dimer saddle dynamics for finding saddle point and solution
  landscape}, Japan Journal of Industrial and Applied Mathematics, 40 (2023),
  pp.~1677--1693.

\bibitem{zhang2024probabilistic}
\leavevmode\vrule height 2pt depth -1.6pt width 23pt, {\em Probabilistic error
  estimate for numerical discretization of high-index saddle dynamics with
  inaccurate models}, Annals of Applied Mathematics, 40 (2024), pp.~1--20.

\bibitem{zhang2025cms}
\leavevmode\vrule height 2pt depth -1.6pt width 23pt, {\em Understanding
  high-index saddle dynamics via numerical analysis}, Communications in
  Mathematical Sciences, 23 (2025), pp.~541--560.

\bibitem{Zhang2025Numerical}
{\sc L.~Zhang, X.~Zheng, and S.~Zhu}, {\em Numerical analysis for saddle
  dynamics of some semilinear elliptic problems}.
\newblock arXiv:2508.01534, 2025.

\bibitem{Zhang2024.11.28.625944}
{\sc X.~Zhang, Z.~Li, and L.~Zhang}, {\em Constructing a holistic map of cell
  fate decision by hyper solution landscape}.
\newblock bioRxiv 2024.11.28.625944, 2024.

\bibitem{zhang2025softmatter}
{\sc Y.~Zhang, X.~Gu, Y.~Wang, X.~Xu, and L.~Zhang}, {\em Solution landscape of
  droplets on rough surfaces: wetting transition and directional transport},
  Soft Matter, 21 (2025), pp.~2729--2737.

\bibitem{jcp2023}
{\sc Y.~Zhang, X.~Yang, L.~Zhang, Y.~Li, T.~Zhang, and S.~Sun}, {\em Energy
  landscape analysis for two-phase multi-component nvt flash systems by using
  etd type high-index saddle dynamics}, Journal of Computational Physics, 477
  (2023), p.~111916.

\bibitem{zhou2024nucleation}
{\sc T.~Zhou, L.~Zhang, P.~Zhang, A.-C. Shi, and K.~Jiang}, {\em Nucleation and
  phase transition of decagonal quasicrystals}, The Journal of Chemical
  Physics, 161 (2024), p.~164503.

\end{thebibliography}

\newpage
\appendix % 声明附录开始
\section{Parameter Reference}\label{Parameter Reference}
The following is the detailed parameter list mentioned in Section \ref{parameter}.

\begin{figure}[htbp]
    \centering
    \includegraphics[width=1.0\linewidth]{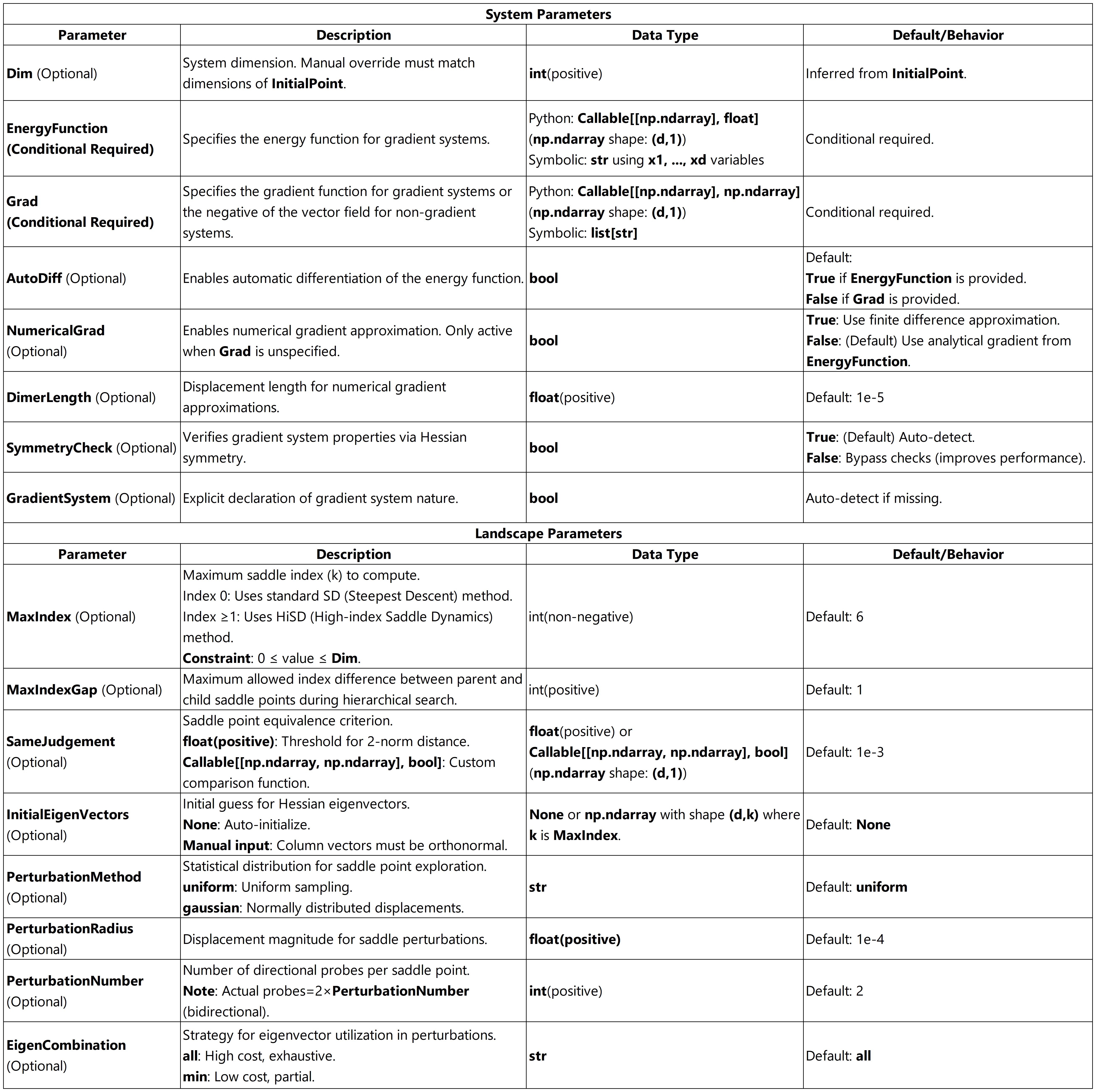}
\end{figure}

\begin{figure}[htbp]
    \centering
    \includegraphics[width=1.0\linewidth]{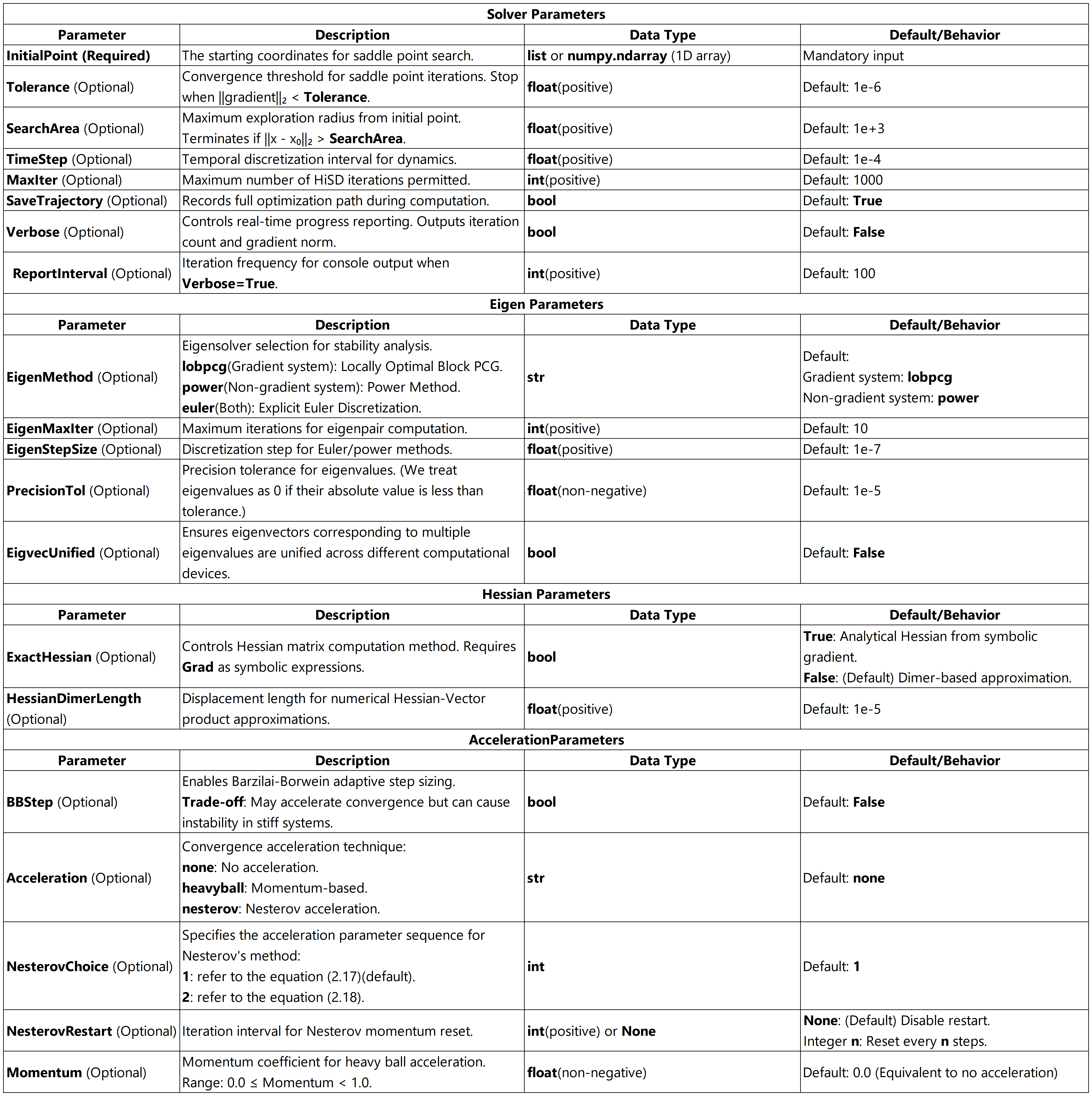}
\end{figure}

\end{document}